\documentclass[a4paper]{amsart}
\usepackage{color}
\usepackage[latin1]{inputenc}
\usepackage[T1]{fontenc}
\usepackage{amsfonts}
\usepackage{amssymb}
\usepackage{amsmath}
\usepackage{amsthm}
\usepackage{graphicx,type1cm,eso-pic,color}
\usepackage{pstricks,pst-plot,pstricks-add}
\usepackage{float}
\newtheorem{theorem}{Theorem}[section]
\newtheorem{lemma}[theorem]{Lemma}
\newtheorem{proposition}[theorem]{Proposition}
\newtheorem{corollary}[theorem]{Corollary}
\newtheorem{definition}[theorem]{Definition}
\newtheorem{assumption}[theorem]{Assumption}
\begin{document}
\setlength\arraycolsep{2pt}
\title{Sensitivities via Rough Paths}
\author{Nicolas MARIE}
\address{Laboratoire Modal'X Universit\'e Paris-Ouest 92000 Nanterre}
\email{nmarie@u-paris10.fr}
\date{}
\maketitle
%


%
\begin{abstract}
Motivated by a problematic coming from mathematical finance, this paper is devoted to existing and additional results of continuity and differentiability of the It\^o map associated to rough differential equations. These regularity results together with Malliavin calculus are applied to sensitivities analysis for stochastic differential equations driven by multidimensional Gaussian processes with continuous paths, especially fractional Brownian motions. Precisely, in that framework, results on computation of greeks for It\^o's stochastic differential equations are extended. An application in mathematical finance, and simulations, are provided.
\end{abstract}
\tableofcontents
\noindent
\textbf{Acknowledgements.} Many thanks to Laure Coutin and Laurent Decreusefond for their advices. That work was supported by ANR Masterie.
%


%
\section{Introduction}
\noindent
Motivated by a problematic coming from mathematical finance, this paper is devoted to existing and additional results of continuity and differentiability of the It\^o map associated to rough differential equations (RDEs). These regularity results together with Malliavin calculus are applied to sensitivities analysis for stochastic differential equations (SDEs) driven by multidimensional Gaussian processes with continuous paths, especially fractional Brownian motions.
\\
\\
In order to state the problematic mentioned above, let's remind some notions of mathematical finance.
\\
\\
Consider a probability space $(\Omega,\mathcal A,\mathbb P)$, a $d$-dimensional Brownian motion $B$ and $\mathbb F := (\mathcal A_t;t\in [0,T])$ the natural filtration associated to $B$ ($d\in\mathbb N^*$ and $T > 0$).
\\
Consider the financial market consisting of $d + 1$ assets (one risk-free asset and $d$ risky assets) over the filtered probability space $(\Omega,\mathcal A,\mathbb F,\mathbb P)$. At time $t\in [0,T]$, the deterministic price of the risk-free asset is denoted by $S_{t}^{0}$, and prices of the $d$ risky assets are given by the random vector $S_t := (S_{t}^{1},\dots,S_{t}^{d})$.
\\
In a first place, assume that the process $S$ is the solution of a stochastic differential equation, taken in the sense of It\^o :
\begin{displaymath}
S_t =
x +
\int_{0}^{t}\mu\left(S_u\right)du +
\int_{0}^{t}\sigma\left(S_u\right)dB_u
\textrm{ ; }x\in\mathbb R^d
\end{displaymath}
where, $\mu : \mathbb R^d\rightarrow\mathbb R^d$ and $\sigma : \mathbb R^d\rightarrow\mathcal M_d(\mathbb R)$ are (globally) Lipschitz continuous functions.
\\
Let $\mathbb P^*\sim\mathbb P$ be the risk-neutral probability measure of the market (i.e. such that $S^* := S/S^0$ is a $(\mathbb F,\mathbb P^*)$-martingale).
%


%
\begin{theorem}\label{Pricing}
Consider an option of payoff $h\in L^2(\Omega,\mathcal A_T,\mathbb P^*)$. Then, there exists an admissible strategy $\varphi$ such that :
\begin{displaymath}
\forall t\in [0,T]\textrm{$,$ }
V_t(\varphi) =
\mathbb E^*
\left(\left.
\frac{S_{t}^{0}}{S_{0}^{0}}h\right|
\mathcal A_t\right)
\textrm{ $\mathbb P^*$-a.s.}
\end{displaymath}
where $V(\varphi)$ is the associated wealth process.
\end{theorem}
\noindent
Refer to D. Lamberton and B. Lapeyre \cite{LL97}, Theorem 4.3.2 for a proof in the case of the Black-Scholes model.
\\
\\
With notations of Theorem \ref{Pricing}, $V_T(\varphi) = \mathbb E^*(h)$. It is the option's price and, when $h := F(S_T)$ with a function $F$ from $\mathbb R^d$ into $\mathbb R_+$, it is useful to get the existence and an expression of the sensitivities of $V_T(\varphi)$ to perturbations of the initial condition and of the volatility function $\sigma$ in particular :
\begin{displaymath}
\Delta :=
\partial_x\mathbb E^*\left[F(S_{T}^{x})\right]
\textrm{ and }
\mathcal V :=
\partial_{\sigma}
\mathbb E^*\left[F(S_{T}^{\sigma})\right].
\end{displaymath}
In finance, these sensitivities are called greeks. For instance, $\Delta$ involves in a procedure called $\Delta$-hedging that allows to compute the admissible strategy of Theorem \ref{Pricing} (cf. \cite{LL97}, Subsection 4.3.3). However, there are no reasons to use these quantities in finance only. Indeed, they could also be used in pharmacokinetic as mentioned at \cite{MARIE12}, Section 5.
\\
Greeks have been studied by many authors. In \cite{FLLLT99}, E. Fourni\'e et al. have established the existence of greeks and have provided expressions of them via Malliavin calculus by assuming that $\sigma$ satisfied the uniform elliptic condition (cf. Theorem \ref{GRQFournie}). In \cite{GOBET05}, E. Gobet and R. M\"unos have extended these results by assuming that $\sigma$ only satisfied the hypoelliptic condition. On computation of greeks in the Black-Scholes model, refer to P. Malliavin and A. Thalmaier \cite{MT06}, Chapter 2. On sensitivities in models with jumps, refer to N. Privault et al. \cite{PW04} and \cite{LP04}. Finally, via the cubature formula for the Brownian motion, J. Teichmann has provided estimators of Malliavin's weights for computation of greeks (cf. J. Teichmann \cite{TEICH05}). On the regularity of the solution map associated to SDEs taken in the sense of It\^o, refer to H. Kunita \cite{KUN97}.
\\
At the following theorem, $\delta$ is the divergence operator associated to the Brownian motion $B$ (cf. D. Nualart \cite{NUALART06}, Section 1.3).
%


%
\begin{theorem}\label{GRQFournie}
When $b$ and $\sigma$ are differentiable, of bounded and Lipschitz continuous derivatives, and $F\in L^2(\mathbb R^d;\mathbb R_+)$ :
\begin{enumerate}
 \item If $\sigma$ satisfies the uniform elliptic condition (i.e. there exists $\varepsilon > 0$ such that for every $a,b\in\mathbb R^d$, $b^T\sigma^T(a)\sigma(a)b\geqslant\varepsilon\|b\|^2$), then $\Delta$ exists and there exists an adapted $d$-dimensional process $h^{\Delta}$ such that :
 \begin{displaymath}
 \Delta =
 \mathbb E^*\left[F(S_T)\delta(h^{\Delta})\right].
 \end{displaymath}
 \item Let $\tilde\sigma : \mathbb R^d\rightarrow\mathcal M_d(\mathbb R)$ be a function such that, for every $\varepsilon$ belonging to a closed neighborhood of $0$, $\sigma +\varepsilon\tilde\sigma$ satisfies the uniform elliptic condition. Then $\mathcal V$ exists and there exists an (anticipative) $d$-dimensional process $h^{\mathcal V}$ such that :
 \begin{displaymath}
 \mathcal V =
 \mathbb E^*\left[F(S_T)\delta(h^{\mathcal V})\right].
 \end{displaymath}
\end{enumerate}
\end{theorem}
\noindent
Refer to \cite{FLLLT99}, propositions 3.2 and 3.3 for a proof.
\\
\\
Under technical assumptions stated at Subsection 2.3, the main purpose of this paper is to extend Theorem \ref{GRQFournie} for the following SDE, taken in the sense of rough paths introduced by T. Lyons in \cite{LYONS98} :
\begin{displaymath}
X_t =
x +
\int_{0}^{t}
\mu\left(X_s\right)ds +
\int_{0}^{t}
\sigma\left(X_s\right)dW_s
\textrm{ ; }
x\in\mathbb R^d
\end{displaymath}
where, $W$ is a centered $d$-dimensional Gaussian process with continuous paths of finite $p$-variation ($p\geqslant 1$), and $\mu$ and $\sigma$ satisfy the following assumption.
%


%
\begin{assumption}\label{REGcoef}
$\mu$ and $\sigma$ are $[p] + 1$ times differentiable, bounded and of bounded derivatives.
\end{assumption}
\noindent
In order to solve that problematic, subsections 2.1 and 2.2 are devoted to existing and additional results of continuity and differentiability of the It\^o map associated to rough differential equations. In particular, the continuous differentiability of the It\^o map with respect to the collection of vector fields is proved in addition to existing regularity results for the initial condition and for the driving signal, coming from P. Friz and N. Victoir \cite{FV08}, chapters 4 and 11. In order to apply the (probabilistic) integrability results coming from T. Cass, C. Litterer and T. Lyons \cite{CLL11}, tailor-made upper-bounds are provided for each derivative. Subsection 2.3 reminds definitions and results related to the natural geometric rough paths associated to a Gaussian process having a covariance function satisfying the technical Assumption \ref{CSGRPgauss}, called enhanced Gaussian process by P. Friz and N. Victoir. Results of subsections 2.1 and 2.2 are applied together with results coming from \cite{CLL11}, in order to get the (probabilistic) integrability of the solution of a Gaussian RDE and their derivatives. The main problematic is solved at Section 3 by using results of Section 2 together with Malliavin calculus arguments. Simulations of $\Delta$ and $\mathcal V$ are provided at Subsection 4.2.
\\
The fractional Brownian motion (fBm) introduced in \cite{MVN68} by B.B. Mandelbrot and W. Van Ness has been studied by many authors in order to generalize the Brownian motion classically used to modelize risky assets prices process. For instance, the regularity of the process paths and its memory are both continuously controlled via the Hurst parameter $H$ of the fBm. However, the fBm is not a semi-martingale if $H\not= 1/2$ (cf. \cite{NUALART06}, Proposition 5.1.1). In \cite{ROGERS97}, L.C.G. Rogers has shown the existence of arbitrages for assets having a prices process modelized by a fBm. In order to bypass that difficulty, several approaches have been studied. For instance, in \cite{CHER01}, P. Cheridito assumed that the prices process was modelized by a mixed fractional Brownian motion that is a semi-martingale that involves a fBm. At Subsection 4.1, prices of the risky assets are modelized by a fractional SDE, in which the volatility is modelized by another one. Results of Section 3 are applied in order to show the existence and provide an expression of the sensitivity of the option's price with respect to collection of vector fields of the volatility's equation.
\\
\\
This paper uses many results on rough paths and rough differential equations coming from \cite{FV08} and, T. Lyons and Z. Qian \cite{LYONS02}. It also involves Malliavin calculus results coming from \cite{NUALART06}.
\\
Notations and, short definitions and results, used throughout the paper are stated below. However, in a sake of rigor, original well stated results of the literature are cited throughout the paper.
\\
\\
\textbf{Notations (general) :}
\begin{itemize}
 \item $\mathbb R^e$ and $\mathbb R^d$ ($e,d\in\mathbb N^*$) are equipped with their euclidean norms, both denoted by $\|.\|$.
 \item The canonical basis of $\mathbb R^d$ is denoted by $(e_1,\dots,e_d)$. With respect to that basis, for $k = 1,\dots,d$ ; the $k$-th component of any vector $u\in\mathbb R^d$ is denoted by $u^k$.
 \item The closed ball of $\mathbb R^d$ for $\|.\|$, of center $a\in\mathbb R^d$ and of radius $r > 0$, is denoted by $B(a,r)$.
 \item The usual matrix (resp. operator) norm on $\mathcal M_{e,d}(\mathbb R)$ (resp. $\mathcal L(\mathbb R^e;\mathbb R^d)$) is denoted by $\|.\|_{\mathcal M}$ (resp. $\|.\|_{\mathcal L}$).
 \item Consider $0\leqslant s < t\leqslant T$. The set of dissections of $[s,t]$ is denoted by $D_{s,t}$. In particular, $D_T := D_{0,T}$.
 \item $\Delta_T := \{(u,v)\in\mathbb R^2 : 0\leqslant u < v\leqslant T\}$.
 \item The space of continuous (resp. continuously differentiable) functions from $[s,t]$ into $\mathbb R^d$ is denoted by $C^0([s,t];\mathbb R^d)$ (resp. $C^1([s,t];\mathbb R^d)$), and is equipped with the uniform norm $\|.\|_{\infty;s,t}$.
 \item Differentiability means differentiability in the sense of Fr\'echet in general (cf. H. Cartan \cite{CART07}, Chapter I.2).
 \item Consider two Banach spaces $E$ and $F$. Let $\varphi : E\rightarrow F$ be a derivable map at point $x\in E$, in the direction $h\in E$ ; the associated directional derivative is denoted by :
 \begin{displaymath}
 D_h\varphi(x) :=
 \lim_{\varepsilon\rightarrow 0}
 \frac{\varphi(x +\varepsilon h) -\varphi(x)}{\varepsilon}
 \textrm{ in }F.
 \end{displaymath}
\end{itemize}
\textbf{Notations (rough paths) :}
\begin{itemize}
 \item Consider $p\geqslant 1$ and $\alpha\in ]0,1]$. The space of continuous functions of finite $p$-variation (resp. $\alpha$-H\"older continuous functions) from $[s,t]$ into $\mathbb R^d$ is denoted by
 \begin{displaymath}
 C^{p\textrm{-var}}\left([s,t];\mathbb R^d\right) :=
 \left\{y\in C^0\left([s,t];\mathbb R^d\right) :
 \sup_{
 D = \left\{r_k\right\}\in D_{s,t}}
 \sum_{k = 1}^{|D| - 1}
 \|y_{r_{k + 1}} - y_{r_k}\|^p < \infty
 \right\}
 \end{displaymath}
 (resp. $C^{\alpha\textrm{-h\"ol}}([s,t];\mathbb R^d)$, that is a subset of $C^{1/\alpha\textrm{-var}}([s,t];\mathbb R^d)$), and is equipped with the $p$-variation distance $d_{p\textrm{-var};s,t}$ (resp. the $\alpha$-H\"older distance $d_{\alpha\textrm{-h\"ol};s,t}$). Refer to \cite{FV08}, chapters 5 and 8 about these spaces.
 \item Consider $N\in\mathbb N^*$ and $y : [0,T]\rightarrow\mathbb R^d$ a continuous function of finite $1$-variation. The step-$N$ tensor algebra over $\mathbb R^d$ is denoted by
 \begin{displaymath}
 T^N(\mathbb R^d) :=
 \mathbb R\oplus\mathbb R^d\oplus\dots\oplus(\mathbb R^d)^{\otimes N},
 \end{displaymath}
 the step-$N$ signature of $y$ is denoted by
 \begin{displaymath}
 S_N(y) :=
 \left(1,\int_{0}^{.}dy_r,\dots,\int_{0 < r_1 < \dots < r_N < .}dy_{r_1}\otimes\dots\otimes dy_{r_N}\right),
 \end{displaymath}
 and the step-$N$ free nilpotent group over $\mathbb R^d$ is denoted by
 \begin{displaymath}
 G^N(\mathbb R^d) :=
 \left\{S_N(y)_1 ; y\in C^{1\textrm{-var}}([0,1];\mathbb R^d)\right\}.
 \end{displaymath}
 Refer to \cite{FV08}, Chapter 7.
 \item For $k = 0,\dots,N$, the $(k + 1)$-th component of any $X\in T^N(\mathbb R^d)$ is denoted by $X^k$.
 \item The space of geometric $p$-rough paths is denoted by
 \begin{displaymath}
 G\Omega_{p,T}(\mathbb R^d) :=
 \overline{
 \left\{
 S_{[p]}(y) ; y\in C^{1\textrm{-var}}([0,T];\mathbb R^d)
 \right\}}^{d_{p\textrm{-var};T}},
 \end{displaymath}
 and is equipped with the $p$-variation distance $d_{p\textrm{-var};T}$, or with the uniform distance $d_{\infty;T}$, associated to the Carnot-Carath\'eodory's distance. Refer to \cite{FV08}, Chapter 9.
 \item The closed ball of $G\Omega_{p,T}(\mathbb R^d)$ for $d_{p\textrm{-var};T}$, of center $Y\in G\Omega_{p,T}(\mathbb R^d)$ and of radius $r > 0$, is denoted by $B_{p,T}(Y,r)$.
 \item For every $Y\in G\Omega_{p,T}(\mathbb R^d)$, $\omega_{Y,p} : (s,t)\in\bar\Delta_T\mapsto\|Y\|_{p\textrm{-var};s,t}^{p}$ is a control. Refer to \cite{FV08}, Chapter 1 about properties of controls.
 \item Consider $q\geqslant 1$ such that $1/p + 1/q > 1$, $Y\in G\Omega_{p,T}(\mathbb R^d)$ and $h\in G\Omega_{q,T}(\mathbb R^e)$. The geometric $p$-rough path over $(Y^1,h^1)$ provided at \cite{FV08}, Theorem 9.26 is denoted by $S_{[p]}(Y\oplus h)$. The translation of $Y$ by $h$ provided at \cite{FV08}, Theorem 9.34 is denoted by $T_hY$.
 \item Consider $\gamma > 0$. The space of collections of $\gamma$-Lipschitz (resp. locally $\gamma$-Lipschitz) vector fields on $\mathbb R^e$ is denoted by $\textrm{Lip}^{\gamma}(\mathbb R^e;\mathbb R^d)$ (resp. $\textrm{Lip}_{\textrm{loc}}^{\gamma}(\mathbb R^e;\mathbb R^d)$) (cf. \cite{FV08}, Definition 10.2). $\textrm{Lip}^{\gamma}(\mathbb R^e;\mathbb R^d)$ is equipped with the $\gamma$-Lipschitz norm $\|.\|_{\textrm{lip}^{\gamma}}$ such that, for every $V\in\textrm{Lip}^{\gamma}(\mathbb R^e;\mathbb R^d)$,
 \begin{displaymath}
 \|V\|_{\textrm{lip}^{\gamma}} :=
 \max\left\{
 \|V\|_{\infty},
 \|DV\|_{\infty},\dots,
 \|D^{\lfloor\gamma\rfloor}V\|_{\infty},
 \|D^{\lfloor\gamma\rfloor}V\|_{\{\gamma\}\textrm{-h\"ol}}
 \right\}.
 \end{displaymath}
 \item The closed ball of $\textrm{Lip}^{\gamma}(\mathbb R^e;\mathbb R^d)$ for $\|.\|_{\textrm{lip}^{\gamma}}$, of center $V\in\textrm{Lip}^{\gamma}(\mathbb R^e;\mathbb R^d)$ and of radius $r > 0$, is denoted by $B_{\textrm{Lip}^{\gamma}}(V,r)$.
 \item Consider $\varepsilon > 0$, $I$ a compact interval of $[0,T]$, a control $\omega :\bar\Delta_T\rightarrow\mathbb R_+$ and $Y\in G\Omega_{p,T}(\mathbb R^d)$. Let's put :
 \begin{displaymath}
 M_{\varepsilon,I,\omega} :=
 \sup_{
 \begin{tiny}
 \begin{array}{rcl}
  D = \{r_k\} & \in & D_I\\
  \omega\left(r_k,r_{k + 1}\right) & \leqslant & \varepsilon
 \end{array}
 \end{tiny}
 }
 \sum_{k = 1}^{|D| - 1}
 \omega\left(r_k,r_{k + 1}\right),
 \end{displaymath}
 $M_{\varepsilon,I,p}(Y) := M_{\varepsilon,I,\omega_{Y,p}}$ and
 \begin{displaymath}
 N_{\varepsilon,I,p}(Y) :=
 \sup\left\{n\in\mathbb N : \tau_n\leqslant\sup(I)\right\}
 \end{displaymath}
 where, $\tau_0 :=\inf(I)$ and for every $n\in\mathbb N$,
 \begin{displaymath}
 \tau_{n + 1} :=
 \inf\left\{t\in I :
 \|Y\|_{p\textrm{-var};\tau_n,t}^{p}\geqslant\varepsilon
 \textrm{ and }
 t > \tau_n\right\}
 \wedge
 \sup(I).
 \end{displaymath}
 In the sequel, $I := [0,T]$.
 \item Consider $\gamma > p$ and $V\in\textrm{Lip}_{\textrm{loc}}^{\gamma}(\mathbb R^e;\mathbb R^d)$ satisfying the $p$-non explosion condition (i.e. $V$ and $D^{[p]}V$ are respectively globally Lipschitz continuous and $(\gamma - [p])$-H\"older continuous on $\mathbb R^e$). The unique solution of $dX = V(X)d\mathbb W$ with initial condition $X_0\in G^{[p]}(\mathbb R^e)$ or $X_0\in\mathbb R^e$, is denoted by $\pi_V(0,X_0;\mathbb W)$.
 \item By \cite{FV08}, Exercice 10.55, if $V$ is a collection of affine vector fields and $\omega :\bar\Delta_T\rightarrow\mathbb R_+$ is a control satisfying $\|\mathbb W\|_{p\textrm{-var};s,t}\leqslant\omega^{1/p}(s,t)$ for every $(s,t)\in\bar\Delta_T$, there exists a constant $C_1 > 0$, not depending on $X_0\in\mathbb R^e$ and $\mathbb W$, such that :
 \begin{displaymath}
 \left\|
 \pi_V(0,x_0;\mathbb W)\right\|_{\infty;T}
 \leqslant
 C_1(1 + \|x_0\|)e^{C_1M_{1,I,\omega}}.
 \end{displaymath}
 By \cite{FV08}, Theorem 10.36, if $V\in\textrm{Lip}^{\gamma}(\mathbb R^e;\mathbb R^d)$, there exists a constant $C_2 > 0$, not depending on $X_0\in G^{[p]}(\mathbb R^e)$, $V$ and $\mathbb W$, such that for every $(s,t)\in\bar\Delta_T$,
 \begin{displaymath}
 \|\pi_V(0,X_0;\mathbb W)\|_{p\textrm{-var};s,t}
 \leqslant
 C_2\left(\|V\|_{\textrm{lip}^{\gamma - 1}}\|\mathbb W\|_{p\textrm{-var};s,t}\vee
 \|V\|_{\textrm{lip}^{\gamma - 1}}^{p}\|\mathbb W\|_{p\textrm{-var};s,t}^{p}\right).
 \end{displaymath}
 By \cite{FV08}, Theorem 10.47, if $V\in\textrm{Lip}^{\gamma}(\mathbb R^e;\mathbb R^d)$, there exists a constant $C_3 > 0$, not depending on $V$ and $\mathbb W$, such that for every $(s,t)\in\bar\Delta_T$,
 \begin{displaymath}
 \left\|
 \int V(\mathbb W)d\mathbb W
 \right\|_{p\textrm{-var};s,t}
 \leqslant
 C_3\|V\|_{\textrm{lip}^{\gamma - 1}}
 \left(\|\mathbb W\|_{p\textrm{-var};s,t}\vee
 \|\mathbb W\|_{p\textrm{-var};s,t}^{p}\right).
 \end{displaymath}
\end{itemize}
\textbf{Notations (Gaussian stochastic analysis) :}
\begin{itemize}
 \item For every $t\in [0,T]$, $[0,t]$ is equipped with the Borel $\sigma$-algebra $\mathcal B_t$ associated to the usual topology on $[0,t]$.
 \item $\mathbb R^d$ is equipped with the Borel $\sigma$-algebra associated to the usual euclidean topology on $\mathbb R^d$, and $G^{[p]}(\mathbb R^d)$ is equipped with the Borel $\sigma$-algebra associated to the Carnot-Carath\'eodory's topology on $G^{[p]}(\mathbb R^d)$. These $\sigma$-algebras are both denoted by $\mathcal B$.
 \item Let $W$ be a $d$-dimensional centered Gaussian process with continuous paths. Its Cameron-Martin's space is denoted by
 \begin{displaymath}
 H^1 :=
 \left\{
 h\in C^0([0,T];\mathbb R) :
 \exists Z\in\mathcal W
 \textrm{ s.t. }
 \forall t\in [0,T]\textrm{, }
 h_t =
 \mathbb E(W_tZ)
 \right\}
 \end{displaymath}
 with
 \begin{displaymath}
 \mathcal W :=
 \overline{\textrm{span}
 \left\{W_t,
 t\in [0,T]\right\}}^{L^2}
 \end{displaymath}
 (cf. \cite{FV08}, Subsection 15.2.2 and Section 15.3), its reproducing kernel Hilbert space is denoted by $H$, and the Wiener integral with respect to $W$ defined on $H$ is denoted by $\mathbf W$ (cf. \cite{NUALART06}, Section 1.1).
 \item The Malliavin derivative associated to $\mathbf W$ is denoted by $\mathbf D$ for the $\mathbb R^d$-valued (resp. $H$-valued) random variables, and its domain in $L^2(\Omega)$ (resp. $L^2(\Omega;H)$) is denoted by $\mathbb D^{1,2}$ (resp. $\mathbb D^{1,2}(H)$) (cf. \cite{NUALART06}, Section 1.2).
 \item For $\mathbb R^d$-valued random variables, the divergence associated to $\mathbf D$ is denoted by $\delta$, and its domain in $L^2(\Omega;H)$ is denoted by $\textrm{dom}(\delta)$ (cf. \cite{NUALART06}, Section 1.3).
 \item The vector space of $\mathbb R^d$-valued (resp. $H$-valued) random variables locally derivable in the sense of Malliavin is denoted by $\mathbb D_{\textrm{loc}}^{1,2}$ (resp. $\mathbb D_{\textrm{loc}}^{1,2}(H)$) (cf. \cite{NUALART06}, Subsection 1.3.5).
\end{itemize}
%


%
\section{Regularity of the It\^o map : existing and additional results}
\noindent
This section is devoted to the regularity of the It\^o map associated to a RDE. On one hand, results of continuity and differentiability of the It\^o map with respect to the initial condition and the driving signal coming from \cite{FV08}, Chapter 11 are stated. In addition, the continuous differentiability of the It\^o map with respect to the collection of vector fields is proved. On the other hand, in order to apply the integrability results coming from \cite{CLL11}, tailor-made upper-bounds are provided for each derivative.
\\
\\
Let's first synthesize existing continuity results for the It\^o map and the rough integral.
%


%
\begin{theorem}\label{ITORDEcontinuite}
Consider $R > 0$ :
\begin{enumerate}
 \item The It\^o map $(X_0,\mathbb W,V)\mapsto\pi_V(0,X_0;\mathbb W)$ is uniformly continuous from
 \begin{displaymath}
 G^{[p]}(\mathbb R^e)\times
 B_{p,T}(\mathbf 1,R)\times
 \textrm{Lip}^{\gamma}(\mathbb R^e;\mathbb R^d)
 \textrm{ into }
 G\Omega_{p,T}(\mathbb R^d).
 \end{displaymath}
 \item The map
 \begin{displaymath}
 \mathbb J : (\mathbb W,V)\longmapsto
 \int V(\mathbb W)d\mathbb W
 \end{displaymath}
 is uniformly continuous from
 \begin{displaymath}
 B_{p,T}(\mathbf 1,R)\times
 \textrm{Lip}^{\gamma - 1}(\mathbb R^d;\mathbb R^d)
 \textrm{ into }
 G\Omega_{p,T}(\mathbb R^d).
 \end{displaymath}
\end{enumerate}
In both cases, the uniform continuity holds when $B_{p,T}(\mathbf 1,R)$ and $G\Omega_{p,T}(\mathbb R^d)$ are equipped with the uniform distance $d_{\infty;T}$.
\end{theorem}
\noindent
Refer to \cite{FV08}, corollaries 10.39,40,48 for a proof.
\\
\\
\textbf{Remark.} Consider $x_0\in\mathbb R^e$, $\mathbb W\in G\Omega_{p,T}(\mathbb R^d)$ and $V := (V_1,\dots,V_d)$ a collection of affine vector fields on $\mathbb R^e$. By \cite{FV08}, Theorem 10.53, $\pi_V(0,x_0;\mathbb W)_t$ belongs to the ball $B(0;R(x_0,\mathbb W))$ of $\mathbb R^e$ for every $t\in [0,T]$, where
\begin{displaymath}
R(x_0,\mathbb W) :=
C(1 + \|x_0\|)e^{C\|\mathbb W\|_{p\textrm{-var};T}^{p}}
\end{displaymath}
and $C > 0$ is a constant not depending on $x_0$ and $\mathbb W$. Moreover, for every $\tilde x_0\in\mathbb R^e$ and every $\widetilde{\mathbb W}\in G\Omega_{p,T}(\mathbb R^d)$,
\begin{displaymath}
\|\tilde x_0\|\leqslant
\|x_0\|
\textrm{ and }
\|\widetilde{\mathbb W}\|_{p\textrm{-var};T}\leqslant
\|\mathbb W\|_{p\textrm{-var};T}
\Longrightarrow
R(\tilde x_0,\widetilde{\mathbb W})\leqslant
R(x_0,\mathbb W).
\end{displaymath}
So, if $\hat V\in\textrm{Lip}^{\gamma}(\mathbb R^e;\mathbb R^d)$ is the collection of vector fields satisfying $\hat V\equiv V$ on $B(0;R(x_0,\mathbb W))$, then
\begin{displaymath}
\pi_V(0,.)\equiv
\pi_{\hat V}(0,.)
\textrm{ on the set }
B(0,\|x_0\|)\times
B_{p,T}(\mathbf 1,\|\mathbb W\|_{p\textrm{-var};T}).
\end{displaymath}
Therefore, by Theorem \ref{ITORDEcontinuite}, the map $\pi_V(0,.)$ is uniformly continuous from
\begin{displaymath}
B(0,\|x_0\|)\times
B_{p,T}(\mathbf 1,\|\mathbb W\|_{p\textrm{-var};T})
\textrm{ into }
C^{p\textrm{-var}}([0,T];\mathbb R^e).
\end{displaymath}
The uniform continuity holds when $B_{p,T}(\mathbf 1,\|\mathbb W\|_{p\textrm{-var};T})$ and $C^{p\textrm{-var}}([0,T];\mathbb R^e)$ are equipped with the uniform distance $d_{\infty;T}$.
\\
\\
The following technical corollary of \cite{FV08}, Theorem 9.26 allows to apply integrability results of \cite{CLL11} to differential equations having a drift term.
%


%
\begin{corollary}\label{MAJcouplage}
Consider $p > q\geqslant 1$ such that $1/p + 1/q > 1$, $Y\in G\Omega_{p,T}(\mathbb R^d)$, $h\in G\Omega_{q,T}(\mathbb R^e)$ and $\varepsilon > 0$. There exists a constant $C > 0$, depending only on $p$ and $q$, such that :
\begin{displaymath}
M_{\varepsilon,I,p}
[S_{[p]}(Y\oplus h)]
\leqslant
C[\|h\|_{q\textrm{-var};T}^{p} +
M_{\varepsilon,I,p}(Y)].
\end{displaymath}
\end{corollary}
%


%
\begin{proof}
On one hand, for every $(s,t)\in\bar\Delta_T$,
\begin{displaymath}
\omega_{Y,p}(s,t) =
\|Y\|_{p\textrm{-var};s,t}
\leqslant
\|S_{[p]}(Y\oplus h)\|_{p\textrm{-var};s,t}.
\end{displaymath}
On the other hand, since $p/q\geqslant 1$ and, $\omega_{Y,p}$ and $\omega_{h,q}$ are two controls :
\begin{displaymath}
\omega =
\|Y\|_{p\textrm{-var}}^{p} +
\|h\|_{q\textrm{-var}}^{p} =
\omega_{Y,p} +
\omega_{h,q}^{p/q}
\end{displaymath}
is also a control.
\\
\\
Then, by \cite{FV08}, Proposition 7.52, there exists a constant $C\geqslant 1$, depending only on $p$ and $q$, such that for every $(s,t)\in\bar\Delta_T$,
\begin{displaymath}
\left\|
S_{[p]}(Y\oplus h)
\right\|_{p\textrm{-var};s,t}^{p}
\leqslant
C\omega(s,t).
\end{displaymath}
In conclusion,
\begin{eqnarray*}
 M_{\varepsilon,I,p}
 \left[S_{[p]}(Y\oplus h)\right] & \leqslant &
 C\sup_{
 \begin{tiny}
 \begin{array}{rcl}
  D = \{r_k\} & \in & D_I\\
  \omega\left(r_k,r_{k + 1}\right) & \leqslant & \varepsilon
 \end{array}
 \end{tiny}
 }
 \sum_{k = 1}^{|D| - 1}
 \omega(r_k,r_{k + 1})\\
 & \leqslant &
 C\left[
 \|h\|_{q\textrm{-var};T}^{p} +
 M_{\varepsilon,I,p}(Y)\right]
\end{eqnarray*}
by super-additivity of the control $\omega_{h,q}^{p/q}$.
\end{proof}
%


%
\subsection{Differentiability of the It\^o map with respect to $x_0$ and $V$}
In order to prove the continuous differentiability of the It\^o map for RDEs with respect to the collection of vector fields, let's first prove it for ODEs.
%


%
\begin{proposition}\label{ITOdiff2}
Consider $\gamma\geqslant 1$, $x_0\in\mathbb R^e$ and a continuous function $w : [0,T]\rightarrow\mathbb R^d$ of finite $1$-variation. Then, the map $V\mapsto\pi_V(0,x_0;w)$ is continuously differentiable from
\begin{displaymath}
\textrm{Lip}^{\gamma}(\mathbb R^e;\mathbb R^d)
\textrm{ into }
C^{1\textrm{-var}}([0,T];\mathbb R^e).
\end{displaymath}
\end{proposition}
%


%
\begin{proof}
In a first step, the derivability of the It\^o map with respect the collection of vector fields is established at every points and in every directions of $\textrm{Lip}^{\gamma}(\mathbb R^d;\mathbb R^e)$. In a second step, via \cite{FV08}, Proposition B.5, the continuous differentiability of that partial It\^o map is proved.
\\
\\
\textbf{Step 1.} Consider $V,\tilde V\in\textrm{Lip}^{\gamma}(\mathbb R^e;\mathbb R^d)$, $\varepsilon\in ]0,1]$, $x^V := \pi_V(0,x_0;w)$ and $y^{V,\tilde V}$ the solution of the following ODE :
\begin{equation}\label{EDOdiff1}
y_{t}^{V,\tilde V} =
\int_{0}^{t}
\langle DV(x_{s}^{V}),y_{s}^{V,\tilde V}\rangle dw_s +
\int_{0}^{t}
\tilde V(x_{s}^{V})dw_s.
\end{equation}
For every $t\in [0,T]$,
\begin{eqnarray*}
 \frac{x_{t}^{V +\varepsilon\tilde V} - x_{t}^{V}}{\varepsilon} -
 y_{t}^{V,\tilde V} & = &
 \int_{0}^{t}\left[
 \frac{V(x_{s}^{V +\varepsilon\tilde V}) - V(x_{s}^{V})}{\varepsilon} -
 \langle DV(x_{s}^{V}),y_{s}^{V,\tilde V}\rangle\right]dw_s +\\
 & &
 \int_{0}^{t}\left[
 \tilde V(x_{s}^{V +\varepsilon\tilde V}) -
 \tilde V(x_{s}^{V})\right]dw_s\\
 & = &
 P_t(\varepsilon) + Q_t(\varepsilon) + R_t(\varepsilon)
\end{eqnarray*}
where,
\begin{eqnarray*}
 P_t(\varepsilon) & := &
 \varepsilon^{-1}\int_{0}^{t}\left[
 V(x_{s}^{V +\varepsilon\tilde V}) - V(x_{s}^{V}) -
 \langle DV(x_{s}^{V}),x_{s}^{V +\varepsilon\tilde V} - x_{s}^{V}\rangle\right]dw_s,\\
 Q_t(\varepsilon) & := &
 \int_{0}^{t}\left[
 \tilde V(x_{s}^{V +\varepsilon\tilde V}) -
 \tilde V(x_{s}^{V})\right]dw_s\textrm{ and}\\
 R_t(\varepsilon) & := &
 \int_{0}^{t}
 \langle DV(x_{s}^{V}),
 \varepsilon^{-1}(x_{s}^{V +\varepsilon\tilde V} - x_{s}^{V}) -
 y_{s}^{V,\tilde V}\rangle dw_s.
\end{eqnarray*}
Firstly, since $V$ is continuously differentiable on $\mathbb R^e$, by \cite{FV08}, Lemma 4.2 :
\begin{eqnarray*}
 \|P_t(\varepsilon)\| & \leqslant &
 \varepsilon^{-1}
 \|w\|_{1\textrm{-var};T}
 \sup_{t\in [0,T]}
 \|V(x_{t}^{V +\varepsilon\tilde V}) - V(x_{t}^V) -
 \langle DV(x_{t}^{V}),x_{t}^{V +\varepsilon\tilde V} - x_{t}^{V}\rangle\|\\
 & \leqslant &
 \eta(\varepsilon)\varepsilon^{-1}
 \|w\|_{1\textrm{-var};T}
 \|x^{V +\varepsilon\tilde V} - x^V\|_{\infty;T}
\end{eqnarray*}
where, $\eta(\varepsilon)\rightarrow 0$ when $\varepsilon\rightarrow 0$.
\\
\\
By \cite{FV08}, Theorem 3.18 :
\begin{equation}\label{MAJdiff1}
\|P(\varepsilon)\|_{\infty;T}
\leqslant
M_3(\varepsilon) :=
2\eta(\varepsilon)e^{3M_1M_2}M_2
\|\tilde V\|_{\infty}
\|w\|_{1\textrm{-var};T}
\end{equation}
with
\begin{displaymath}
M_1 :=
\|V\|_{\textrm{lip}^{\gamma}} +\|\tilde V\|_{\textrm{lip}^{\gamma}}
\geqslant
\|V +\varepsilon\tilde V\|_{\textrm{lip}^1}\vee
\|V\|_{\textrm{lip}^1}
\end{displaymath}
and $M_2 := \|w\|_{1\textrm{-var};T}$.
\\
\\
Secondly, since $\tilde V$ is continuously differentiable and of bounded derivative on $\mathbb R^e$, it is a collection of Lipschitz continuous vector fields. Then, by \cite{FV08}, Theorem 3.18 again :
\begin{equation}\label{MAJdiff2}
\|Q(\varepsilon)\|_{\infty;T}\leqslant
M_4(\varepsilon) :=
2\varepsilon e^{3M_1M_2}M_2
\|\tilde V\|_{\textrm{lip}^{\gamma}}^{2}
\|w\|_{1\textrm{-var};T}.
\end{equation}
Thirdly,
\begin{equation}\label{MAJdiff3}
 \|R_t(\varepsilon)\|\leqslant
 \|V\|_{\textrm{lip}^{\gamma}}
 \int_{0}^{t}
 \left\|
 \frac{x_{s}^{V +\varepsilon\tilde V} - x_{s}^{V}}{\varepsilon} - y_{s}^{V,\tilde V}\right\|
 \|dw_s\|.
\end{equation}
Therefore, inequalities (\ref{MAJdiff1}), (\ref{MAJdiff2}) and (\ref{MAJdiff3}) imply that :
\begin{displaymath}
\left\|
\frac{x_{t}^{V +\varepsilon\tilde V} - x_{t}^{V}}{\varepsilon} - y_{t}^{V,\tilde V}\right\|
\leqslant
M_3(\varepsilon) + M_4(\varepsilon) +
\|V\|_{\textrm{lip}^{\gamma}}
\int_{0}^{t}
\left\|
\frac{x_{s}^{V +\varepsilon\tilde V} - x_{s}^{V}}{\varepsilon} - y_{s}^{V,\tilde V}\right\|
\|dw_s\|.
\end{displaymath}
In conclusion, by Gronwall's lemma :
\begin{eqnarray*}
 \left\|
 \frac{x^{V +\varepsilon\tilde V} - x^V}{\varepsilon} - y^{V,\tilde V}\right\|_{\infty;T}
 & \leqslant &
 \left[M_3(\varepsilon) + M_4(\varepsilon)\right]
 e^{\|V\|_{\textrm{lip}^{\gamma}}\|w\|_{1\textrm{-var};T}}\\
 & &
 \xrightarrow[\varepsilon\rightarrow 0]{} 0.
\end{eqnarray*}
\textbf{Step 2.} The solution of equation (\ref{EDOdiff1}) satisfies :
\begin{displaymath}
D_{\tilde V}x^V =
\pi_A(0,0;.)\circ
\mathbb J[F_{V,\tilde V}(.),.]\circ
(\pi_V(0,x_0;.),.)(w)
\end{displaymath}
where, $A :\mathbb R^e\rightarrow\mathcal L(\mathcal L(\mathbb R^e)\times\mathbb R^e;\mathbb R^e)$ and $F_{V,\tilde V} : \mathbb R^e\times\mathbb R^d\rightarrow\mathcal L(\mathbb R^e\times\mathbb R^d;\mathcal L(\mathbb R^e)\times\mathbb R^e)$ are two collections of vector fields, respectively defined by :
\begin{eqnarray*}
 A(a)(L,b) & := & L.a + b\textrm{ and}\\
 F_{V,\tilde V}(a,a')(b,b') & := &
 (\langle DV(a),.\rangle b' ; \tilde V(a) b')
\end{eqnarray*}
for every $a,b\in\mathbb R^e$, $a',b'\in\mathbb R^d$ and $L\in\mathcal L(\mathbb R^e)$.
\\
\\
Firstly, by the second point of Theorem \ref{ITORDEcontinuite}, the map $\mathbb J$ is uniformly continuous on every bounded sets of
\begin{displaymath}
C^{1\textrm{-var}}\left([0,T];\mathbb R^e\times\mathbb R^d\right)\times
C^{1\textrm{-var}}\left([0,T];\mathbb R^e\times\mathbb R^d\right).
\end{displaymath}
Secondly, the map $(V,\tilde V,a)\mapsto F_{V,\tilde V}(a)$ is uniformly continuous on every bounded sets of
\begin{displaymath}
\textrm{Lip}^{\gamma}\left(\mathbb R^e;\mathbb R^d\right)\times
\textrm{Lip}^{\gamma}\left(\mathbb R^e;\mathbb R^d\right)\times
\mathbb R^e\times\mathbb R^d
\end{displaymath}
by construction.
\\
\\
Thirdly, maps $\pi_A(0,0;.)$ and $V\mapsto\pi_V(0,x_0;w)$ are respectively uniformly continuous on every bounded sets of
\begin{displaymath}
C^{1\textrm{-var}}\left([0,T];\mathcal L(\mathbb R^e)\times\mathbb R^e\right)\textrm{ and }
\textrm{Lip}^{\gamma}\left(\mathbb R^e;\mathbb R^d\right)
\end{displaymath}
by Theorem \ref{ITORDEcontinuite} and its remark.
\\
\\
Therefore, by composition, the map $(V,\tilde V)\mapsto D_{\tilde V}x^V$ is uniformly continuous on every bounded sets of
\begin{displaymath}
\textrm{Lip}^{\gamma}\left(\mathbb R^e;\mathbb R^d\right)\times
\textrm{Lip}^{\gamma}\left(\mathbb R^e;\mathbb R^d\right).
\end{displaymath}
In conclusion, by \cite{FV08}, Proposition B.5, the map $V\mapsto\pi_V(0,x_0;w)$ is continuously differentiable from
\begin{displaymath}
\textrm{Lip}^{\gamma}(\mathbb R^e;\mathbb R^d)
\textrm{ into }
C^{1\textrm{-var}}([0,T];\mathbb R^e).
\end{displaymath}
\end{proof}
%


%
\begin{theorem}\label{CIVFdiff}
Consider $\mathbb W\in G\Omega_{p,T}(\mathbb R^d)$ :
\begin{enumerate}
 \item Let $V := (V_1,\dots,V_d)$ be a collection of $\gamma$-Lipschtiz vector fields on $\mathbb R^e$. The map $x_0\mapsto\pi_V(0,x_0;\mathbb W)$ is continuously differentiable from
 \begin{displaymath}
 \mathbb R^e
 \textrm{ into }
 C^{p\textrm{-var}}\left([0,T];\mathbb R^e\right).
 \end{displaymath}
 For every $t\in [0,T]$, the Jacobian matrix of $\pi_V(0,.;\mathbb W)_t$ at point $x_0\in\mathbb R^e$ is denoted by $J_{t\leftarrow 0}^{x_0,\mathbb W}$.
 \\
 \\
 Moreover, for every $\varepsilon > 0$, there exists a constant $C_1 > 0$ only depending on $p$, $\gamma$, $\varepsilon$ and $\|V\|_{\textrm{lip}^{\gamma}}$, such that for every $x_0\in\mathbb R^e$,
 \begin{displaymath}
 \|J_{.\leftarrow 0}^{x_0,\mathbb W}\|_{\infty;T}\leqslant
 C_1e^{C_1M_{\varepsilon,I,p}(\mathbb W)}.
 \end{displaymath}
 \item For every $t\in [0,T]$, $J_{t\leftarrow 0}^{x_0,\mathbb W}$ is an invertible matrix. Moreover, for every $\varepsilon > 0$, there exists a constant $C_2 > 0$ only depending on $p$, $\gamma$, $\varepsilon$ and $\|V\|_{\textrm{lip}^{\gamma}}$, such that for every $x_0\in\mathbb R^e$,
 \begin{displaymath}
 \|(J_{.\leftarrow 0}^{x_0,\mathbb W})^{-1}\|_{\infty;T}\leqslant
 C_2e^{C_2M_{\varepsilon,I,p}(\mathbb W)}.
 \end{displaymath}
 \item Consider $x_0\in\mathbb R^e$. The map $V\mapsto\pi_V(0,x_0;\mathbb W)$ is continuously differentiable from
 \begin{displaymath}
 \textrm{Lip}^{\gamma}(\mathbb R^e;\mathbb R^d)
 \textrm{ into }
 C^{p\textrm{-var}}([0,T];\mathbb R^e).
 \end{displaymath}
 Moreover, for every $R > 0$ and $V,\tilde V\in B_{\textrm{Lip}^{\gamma}}(0,R)$, there exists two constants $\eta > 0$ and $C_3 > 0$, depending on $R$ (continuously) and not on $\mathbb W$, such that :
 \begin{displaymath}
 \|\partial_V\pi_V(0,x_0;\mathbb W).\tilde V\|_{\infty;T}
 \leqslant
 C_3e^{C_3M_{\eta,I,p}(\mathbb W)}.
 \end{displaymath}
\end{enumerate}
\end{theorem}
%


%
\begin{proof}
Refer to \cite{FV08}, theorems 11.3-6 for a proof of the It\^o map's continuous differentiability with respect to the initial condition, and refer to \cite{CLL11}, Corollary 3.4 for the upper-bound provided at the first point for $\|J_{.\leftarrow 0}^{x_0,\mathbb W}\|_{\infty;T}$ ; $x_0\in\mathbb R^e$.
\\
\\
Let $I$ be the identity matrix of $\mathcal M_e(\mathbb R)$. Proofs of points 1 and 2 are similar because if $w : [0,T]\rightarrow\mathbb R^d$ is a continuous function of finite $1$-variation, then
\begin{eqnarray*}
 J_{t\leftarrow 0}^{x_0,w} & = & I + \int_{0}^{t}\langle DV[\pi_V(0,x_0;\mathbb W)_s],J_{s\leftarrow 0}^{x_0,w}\rangle dw_s\textrm{ and}\\
 (J_{t\leftarrow 0}^{x_0,w})^{-1} & = & I - \int_{0}^{t}\langle DV[\pi_V(0,x_0;\mathbb W)_s],(J_{s\leftarrow 0}^{x_0,w})^{-1}\rangle dw_s
\end{eqnarray*}
as mentioned at the proof of \cite{FV08}, Proposition 4.11.
\\
\\
The proof of the third point is detailed. In a first step, the continuous differentiability of the It\^o map with respect to the collection of vector fields is established. In a second step, in order to apply integrability results coming from \cite{CLL11}, a tailor-maid upper-bound for the derivative of the It\^o map with respect to $V$ is provided.
\\
\\
\textbf{Step 1.} Since $\mathbb W\in G\Omega_{p,T}(\mathbb R^d)$, there exists a sequence $(w^n,n\in\mathbb N)$ of functions belonging to $C^{1\textrm{-var}}([0,T];\mathbb R^d)$ and satisfying :
\begin{equation}\label{DVFCVgrp}
\lim_{n\rightarrow\infty}
d_{p\textrm{-var};T}\left[
S_{[p]}(w^n)_{0,.},\mathbb W\right] = 0.
\end{equation}
Consider $n\in\mathbb N$, $\mathbb W^n := S_{[p]}(w^n)_{0,.}$, $x_0\in\mathbb R^e$, $a := (x_0,0)$,
\begin{displaymath}
X_0 :=
\left(1,a,\dots,\frac{a^{\otimes [p]}}{[p]!}\right)\in T^{[p]}\left(\mathbb R^{e + 1}\right)
\end{displaymath}
and $V,\tilde V\in\textrm{Lip}^{\gamma}(\mathbb R^e;\mathbb R^d)$.
\\
\\
By Proposition \ref{ITOdiff2}, the map $\pi_.(0,x_0;w^n)$ is continuously differentiable from
\begin{displaymath}
\textrm{Lip}^{\gamma}(\mathbb R^e;\mathbb R^d)
\textrm{ into }
C^{1\textrm{-var}}([0,T];\mathbb R^e).
\end{displaymath}
In particular, $\partial_V\pi_V(0,x_0;w^n).\tilde V = \varphi(\mathbb W^n,V,\tilde V)$ with
\begin{displaymath}
\varphi(.,V,\tilde V) :=
\pi_A(0,0;.)\circ\mathbb J(. ; F_{V,\tilde V})\circ\pi_{F_V}(0,X_0;.)
\end{displaymath}
where,
\begin{eqnarray*}
 A : \mathbb R^e
 & \longrightarrow &
 \mathcal L(\mathcal L(\mathbb R^e)\times\mathbb R^e ; \mathbb R^e),\\
 F_{V,\tilde V} : \mathbb R^e\times\mathbb R^d
 & \longrightarrow &
 \mathcal L(\mathbb R^e\times\mathbb R^d ; \mathcal L(\mathbb R^e)\times\mathbb R^e)
 \textrm{ and}\\
 F_V : \mathbb R^e
 & \longrightarrow &
 \mathcal L(\mathbb R^d ; \mathbb R^e\times\mathbb R^d)
\end{eqnarray*}
are three collections of vector fields, respectively defined by :
\begin{eqnarray*}
 A(a)(L,b) & := & L.a + b,\\
 F_{V,\tilde V}(a,a')(b,b') & := & (\langle DV(a),.\rangle b' ; \tilde V(a)b')\textrm{ and}\\
 F_V(a)b' & := & (V(a)b',b')
\end{eqnarray*}
for every $a,b\in\mathbb R^e$, $a',b'\in\mathbb R^d$ and $L\in\mathcal L(\mathbb R^e)$.
\\
\\
Consider $\varepsilon\in ]0,1]$. By Taylor's formula applied to $\pi_.(0,x_0;\mathbb W^n)$ between $V$ and $V +\varepsilon\tilde V$, and \cite{FV08}, Definition 10.17 (RDE's solution(s)) :
\begin{equation}\label{DVFexp}
\pi_{V +\varepsilon\tilde V}(0,x_0;\mathbb W) -
\pi_V(0,x_0;\mathbb W) =
\lim_{n\rightarrow\infty}
\int_{0}^{\varepsilon}
\varphi(\mathbb W^n,V +\theta\tilde V,\tilde V)d\theta
\end{equation}
uniformly.
\\
\\
Via Lebesgue's theorem and \cite{FV08}, Proposition B.1, let show that the derivative of $\pi_.(0,x_0;\mathbb W)$ at point $V$, in the direction $\tilde V$, exists in $C^{p\textrm{-var}}([0,T];\mathbb R^e)$ equipped with the norm $\|.\|_{p\textrm{-var};T}$ and coincides with $\varphi(\mathbb W,V,\tilde V)$.
\\
\\
On one hand, the continuity results of Theorem \ref{ITORDEcontinuite} imply that :
\begin{displaymath}
\forall\theta\in ]0,1]\textrm{, }
\varphi(\mathbb W^n,V +\theta\tilde V,\tilde V)
\xrightarrow[n\rightarrow\infty]{}
\varphi(\mathbb W,V +\theta\tilde V,\tilde V)
\end{displaymath}
in $C^{p\textrm{-var}}([0,T];\mathbb R^e)$ equipped with $\|.\|_{\infty;T}$.
\\
\\
On the other hand, by applying successively \cite{FV08}, theorems 10.47 and 10.36, for every $\theta\in ]0,1]$ and every $(s,t)\in\bar\Delta_T$,
\begin{eqnarray*}
 \omega_{1}^{1/p}(s,t;n;\theta) & := &
 \left\|\int F_{V +\theta\tilde V,\tilde V}
 \left[\pi_{F_{V +\theta\tilde V}}(0,X_0;\mathbb W^n)\right]
 d\pi_{F_{V +\theta\tilde V}}(0,X_0;\mathbb W^n)\right\|_{p\textrm{-var};s,t}\\
 & \leqslant &
 \omega_{2}^{1/p}(s,t;n)
\end{eqnarray*}
with
\begin{displaymath}
\omega_{2}^{1/p}(s,t;n) :=
\omega_{3}^{1/p}(s,t;n)\vee
\omega_3(s,t;n)\vee
\omega_{3}^{p}(s,t;n)
\end{displaymath}
and
\begin{displaymath}
\omega_3(s,t;n) :=
\eta_1
\|\mathbb W^n\|_{p\textrm{-var};s,t}^{p}
\end{displaymath}
where $\eta_1 > 0$ is depending on $V$ and $\tilde V$, but not on $\mathbb W^n$ and $\theta$.
\\
\\
By \cite{FV08}, Exercice 10.55, there exists a constant $C_4 > 0$, not depending on $\mathbb W^n$ and $\theta$, such that :
\begin{eqnarray*}
 \left\|\varphi(\mathbb W^n,V +\theta\tilde V,\tilde V)\right\|_{\infty;T} & \leqslant &
 C_4\exp\left[C_4
 \sup_{
 \begin{tiny}
 \begin{array}{rcl}
  D = \{r_k\} & \in & D_I\\
  \omega_2\left(r_k,r_{k + 1};n\right) & \leqslant & 1
 \end{array}
 \end{tiny}
 }
 \sum_{k = 1}^{|D| - 1}
 \omega_2(r_k,r_{k + 1};n)\right]\\
 & = &
 C_4\exp\left[C_4
 \sup_{
 \begin{tiny}
 \begin{array}{rcl}
  D = \{r_k\} & \in & D_I\\
  \omega_3\left(r_k,r_{k + 1};n\right) & \leqslant & 1
 \end{array}
 \end{tiny}
 }
 \sum_{k = 1}^{|D| - 1}
 \omega_3(r_k,r_{k + 1};n)\right],
\end{eqnarray*}
because
\begin{displaymath}
\omega_2(.;n)\equiv\omega_3(.;n)
\textrm{ when }
\omega_2(.;n)\leqslant 1.
\end{displaymath}
By super-additivity of the control $\omega_3(.;n)$ :
\begin{displaymath}
\left\|\varphi(\mathbb W^n,V +\theta\tilde V,\tilde V)\right\|_{\infty;T}
\leqslant
C_4e^{\eta_1C_4\left\|\mathbb W^n\right\|_{p\textrm{-var};T}^{p}}.
\end{displaymath}
In the right-hand side of that inequality, since $\eta_1$ and $C_4$ are not depending on $\mathbb W^n$ and $\theta$, and since
\begin{displaymath}
\sup_{n\in\mathbb N^*}
\left\|\mathbb W^n\right\|_{p\textrm{-var};T}^{p} < \infty
\end{displaymath}
by (\ref{DVFCVgrp}) ;
\begin{displaymath}
\sup_{\theta\in [0,1]}
\sup_{n\in\mathbb N}
\left\|\varphi(\mathbb W^n,V +\theta\tilde V,\tilde V)\right\|_{\infty;T} < \infty
\end{displaymath}
in $C^{p\textrm{-var}}([0,T];\mathbb R^e)$ equipped with $\|.\|_{\infty;T}$.
\\
\\
Therefore, by Lebesgue's theorem and inequality (\ref{DVFexp}) :
\begin{displaymath}
\pi_{V +\varepsilon\tilde V}(0,x_0;\mathbb W) -
\pi_V(0,x_0;\mathbb W) =
\int_{0}^{\varepsilon}
\varphi(\mathbb W,V +\theta\tilde V,\tilde V)d\theta.
\end{displaymath}
Since $\theta\mapsto\varphi(\mathbb W,V +\theta\tilde V,\tilde V)$ is continuous from
\begin{displaymath}
[0,1]
\textrm{ into }
C^{p\textrm{-var}}([0,T];\mathbb R^e)\textrm{ (equipped with $\|.\|_{p\textrm{-var};T}$)}
\end{displaymath}
by Theorem \ref{ITORDEcontinuite} ; by \cite{FV08}, Proposition B.1, the derivative of $\pi_.(0,x_0;\mathbb W)$ at point $V$, in the direction $\tilde V$, exists in $C^{p\textrm{-var}}([0,T];\mathbb R^e)$ equipped with $\|.\|_{p\textrm{-var};T}$ and coincides with $\varphi(\mathbb W,V,\tilde V)$.
\\
\\
Finally, as at the second step of the proof of Proposition \ref{ITOdiff2}, via \cite{FV08}, Proposition B.5 and Lemma 4.2 ; the map $\pi_.(0,x_0;\mathbb W)$ is continuously differentiable from
\begin{displaymath}
\textrm{Lip}^{\gamma}(\mathbb R^e;\mathbb R^d)
\textrm{ into }
C^{p\textrm{-var}}([0,T];\mathbb R^e).
\end{displaymath}
\textbf{Step 2.} Consider $R > 0$ and $V,\tilde V\in B_{\textrm{Lip}^{\gamma}}(0,R)$.
\\
\\
By applying successively \cite{FV08}, theorems 10.47 and 10.36, for every $(s,t)\in\bar\Delta_T$,
\begin{eqnarray*}
 \omega_{4}^{1/p}(s,t) & := &
 \left\|\int F_{V,\tilde V}
 \left[\pi_{F_V}(0,X_0;\mathbb W)\right]
 d\pi_{F_V}(0,X_0;\mathbb W)\right\|_{p\textrm{-var};s,t}\\
 & \leqslant &
 \omega_{5}^{1/p}(s,t)
\end{eqnarray*}
with
\begin{displaymath}
\omega_{5}^{1/p}(s,t) :=
\omega_{6}^{1/p}(s,t)\vee
\omega_6(s,t)\vee
\omega_{6}^{p}(s,t)
\end{displaymath}
and
\begin{displaymath}
\omega_6(s,t) :=
\eta_2
\|\mathbb W\|_{p\textrm{-var};s,t}^{p}
\end{displaymath}
where $\eta_2 > 0$ is depending on $R$ (continuously), but not on $\mathbb W$.
\\
\\
By \cite{FV08}, Exercice 10.55, there exists a constant $C_5 > 0$, not depending on $R$ and $\mathbb W$, such that :
\begin{eqnarray*}
 \left\|\partial_V\pi_V(0,x_0;\mathbb W).\tilde V\right\|_{\infty;T} & \leqslant &
 C_5\exp\left[C_5
 \sup_{
 \begin{tiny}
 \begin{array}{rcl}
  D = \{r_k\} & \in & D_I\\
  \omega_5\left(r_k,r_{k + 1}\right) & \leqslant & 1
 \end{array}
 \end{tiny}
 }
 \sum_{k = 1}^{|D| - 1}
 \omega_5(r_k,r_{k + 1})\right]\\
 & = &
 C_5\exp\left[C_5
 \sup_{
 \begin{tiny}
 \begin{array}{rcl}
  D = \{r_k\} & \in & D_I\\
  \omega_6\left(r_k,r_{k + 1}\right) & \leqslant & 1
 \end{array}
 \end{tiny}
 }
 \sum_{k = 1}^{|D| - 1}
 \omega_6(r_k,r_{k + 1})\right],
\end{eqnarray*}
because
\begin{displaymath}
\omega_5\equiv\omega_6
\textrm{ when }
\omega_5\leqslant 1.
\end{displaymath}
However,
\begin{displaymath}
\sup_{
\begin{tiny}
\begin{array}{rcl}
 D = \{r_k\} & \in & D_I\\
 \omega_6\left(r_k,r_{k + 1}\right) & \leqslant & 1
\end{array}
\end{tiny}
}
\sum_{k = 1}^{|D| - 1}
\omega_6(r_k,r_{k + 1}) =
\eta_2
M_{\eta_{2}^{-1},I,p}(\mathbb W).
\end{displaymath}
Therefore,
\begin{displaymath}
\left\|\partial_V\pi_V(0,x_0;\mathbb W).\tilde V\right\|_{\infty;T}
\leqslant
C_3e^{C_3M_{\eta,I,p}(\mathbb W)}
\end{displaymath}
with $C_3 := C_5(1\vee\eta_2)$ and $\eta := \eta_{2}^{-1}$.
\end{proof}
\noindent
\textbf{Notations.} In the sequel, matrices $J_{t\leftarrow 0}^{x_0,\mathbb W}$ and $(J_{t\leftarrow 0}^{x_0,\mathbb W})^{-1}$ will be respectively denoted by $J_{0\leftarrow t}^{\mathbb W}$ and $J_{t\leftarrow 0}^{\mathbb W}$ in a sake of simplicity. Moreover, for every $(s,t)\in\bar\Delta_T$, let's put :
\begin{displaymath}
J_{s\leftarrow t}^{\mathbb W} :=
J_{s\leftarrow 0}^{\mathbb W}
J_{0\leftarrow t}^{\mathbb W}
\textrm{ and }
J_{t\leftarrow s}^{\mathbb W} :=
J_{t\leftarrow 0}^{\mathbb W}
J_{0\leftarrow s}^{\mathbb W}.
\end{displaymath}
Then,
\begin{displaymath}
J_{s\leftarrow t}^{\mathbb W}
J_{t\leftarrow s}^{\mathbb W} =
J_{t\leftarrow s}^{\mathbb W}
J_{s\leftarrow t}^{\mathbb W} = I.
\end{displaymath}
At the following corollary, upper-bounds provided at the previous theorem are extended to RDEs having a drift term.
%


%
\begin{corollary}\label{MAJsdrift}
Consider $m\in\mathbb N^*$, $p > q\geqslant 1$ such that $1/p + 1/q > 1$, $h : [0,T]\rightarrow\mathbb R^m$ a continuous function of finite $q$-variation, $\mathbb W\in G\Omega_{p,T}(\mathbb R^d)$ and $\mathbb W^h := S_{[p]}(\mathbb W\oplus h)$ :
\begin{enumerate}
 \item Let $V := (V_1,\dots,V_{d + m})$ be a collection of $\gamma$-Lipschitz vector fields on $\mathbb R^e$. For every $\varepsilon > 0$, there exists a constant $C_1 > 0$ depending only on $p$, $q$, $\gamma$, $\varepsilon$ and $\|V\|_{\textrm{lip}^{\gamma}}$, such that for every $x_0\in\mathbb R^e$,
 \begin{displaymath}
 \|J_{.\leftarrow 0}^{\mathbb W^h}\|_{\infty;T}
 \leqslant
 C_1\exp\left[C_1\left[\|h\|_{q\textrm{-var};T}^{p} + M_{\varepsilon,I,p}(\mathbb W)\right]\right].
 \end{displaymath}
 \item Consider $x_0\in\mathbb R^e$. For every $R > 0$ and $V,\tilde V\in B_{\textrm{Lip}^{\gamma}}(0,R)$, there exists two constants $\varepsilon > 0$ and $C_2 > 0$, depending on $R$ but not on $h$ and $\mathbb W$, such that :
 \begin{displaymath}
 \|\partial_V\pi_V(0,x_0;\mathbb W^h).\tilde V\|_{\infty;T}
 \leqslant
 C_2\exp\left[C_2\left[\|h\|_{q\textrm{-var};T}^{p} + M_{\varepsilon,I,p}(\mathbb W)\right]\right].
 \end{displaymath}
\end{enumerate}
\end{corollary}
%


%
\begin{proof}
By Corollary \ref{MAJcouplage}, there exists a constant $C_3 > 0$, depending only on $p$ and $q$, such that for every $\varepsilon > 0$,
\begin{displaymath}
M_{\varepsilon,I,p}(\mathbb W^h)
\leqslant
C_3\left[\|h\|_{q\textrm{-var};T}^{p} + M_{\varepsilon,I,p}(\mathbb W)\right].
\end{displaymath}
Therefore, by Theorem \ref{CIVFdiff} :
\begin{enumerate}
 \item For a given collection of vector fields $V\in\textrm{Lip}^{\gamma}(\mathbb R^e;\mathbb R^{d + m})$ ; for every $\varepsilon > 0$, there exists a constant $C_4 > 0$ depending only on $p$, $\gamma$, $\varepsilon$ and $\|V\|_{\textrm{lip}^{\gamma}}$, such that for every $x_0\in\mathbb R^e$,
 \begin{eqnarray*}
  \|J_{.\leftarrow 0}^{\mathbb W^h}\|_{\infty;T} & \leqslant &
  C_4e^{C_4M_{\varepsilon,I,p}(\mathbb W^h)}\\
  & \leqslant &
  C_1\exp\left[C_1\left[\|h\|_{q\textrm{-var};T}^{p} + M_{\varepsilon,I,p}(\mathbb W)\right]\right]
 \end{eqnarray*}
 with $C_1 := C_4(1\vee C_3)$.
 \item For a given initial condition $x_0\in\mathbb R^e$ ; for every $R > 0$ and $V,\tilde V\in B_{\textrm{Lip}^{\gamma}}(0,R)$, there exists two constants $\varepsilon > 0$ and $C_5 > 0$, depending on $R$ but not on $\mathbb W^h$, such that :
 \begin{eqnarray*}
  \|\partial_V\pi_V(0,x_0;\mathbb W^h).\tilde V\|_{\infty;T}
  & \leqslant &
  C_5e^{C_5M_{\varepsilon,I,p}(\mathbb W^h)}\\
  & \leqslant &
  C_2\exp\left[C_2\left[
  \|h\|_{q\textrm{-var};T}^{p} +
  M_{\varepsilon,I,p}(\mathbb W)\right]\right]
 \end{eqnarray*}
 with $C_2 := C_5(1\vee C_3)$.
\end{enumerate}
\end{proof}
%


%
\subsection{Differentiability of the It\^o map with respect to the driving signal}
Let's first remind the notion of differentiability introduced by P. Friz and N. Victoir on $G\Omega_{p,T}(\mathbb R^d)$.
%


%
\begin{definition}\label{FVdiff}
Consider a Banach space $F$, $p > q\geqslant 1$ such that $1/p + 1/q > 1$, and an open set $U$ of $G\Omega_{p,T}(\mathbb R^d)$. The map $\varphi : G\Omega_{p,T}(\mathbb R^d)\rightarrow F$ is continuously differentiable in the sense of Friz-Victoir on $U$ if and only if, for every $Y\in U$,
\begin{displaymath}
h\in
C^{q\textrm{-var}}([0,T];\mathbb R^d)
\longmapsto
\varphi(T_hY)\in F
\end{displaymath}
is continuously differentiable.
\end{definition}
\noindent
With notations of Definition \ref{FVdiff}, if $\varphi$ is continuously differentiable from $U$ into $F$ in the sense of Friz-Victoir, in particular :
\begin{displaymath}
\forall Y\in U\textrm{, }
\psi^Y : h\in C^{q\textrm{-var}}([0,T];\mathbb R^d)
\longmapsto
\psi^Y(h) = \varphi(T_hY)
\end{displaymath}
is derivable at every points and in every directions of $C^{q\textrm{-var}}([0,T];\mathbb R^d)$.
\\
\\
\textbf{Notation.} For every continuous functions $h : [0,T]\rightarrow\mathbb R^d$ of finite $q$-variation,
\begin{eqnarray*}
 D_{h}^{\textrm{FV}}\varphi(Y) & := &
 D_h\psi^Y(0)\\
 & = &
 \lim_{\varepsilon\rightarrow 0}
 \frac{\varphi(T_{\varepsilon h}Y) -\varphi(T_0Y)}{\varepsilon}.
\end{eqnarray*}
In the sequel, $D^{\textrm{FV}}$ is called the Friz-Victoir (directional) derivative operator.
%


%
\begin{theorem}\label{DFVsignal}
Consider a collection $V := (V_1,\dots,V_d)$ of $\gamma$-Lipschitz vector fields on $\mathbb R^e$ and $x_0\in\mathbb R^e$. The map $\mathbb W\mapsto\pi_V(0,x_0;\mathbb W)$ is continuously differentiable from
\begin{displaymath}
G\Omega_{p,T}(\mathbb R^d)
\textrm{ into }
C^{p\textrm{-var}}([0,T];\mathbb R^e)
\end{displaymath}
in the sense of Friz-Victoir.
\\
\\
Moreover, for every $\mathbb W\in G\Omega_{p,T}(\mathbb R^d)$ and every continous function $h : [0,T]\rightarrow\mathbb R^d$ of finite $q$-variation,
\begin{displaymath}
D_{h}^{FV}\pi_V(0,x_0;\mathbb W) =
\int_{0}^{.}
J_{.\leftarrow s}^{\mathbb W}V\left[
\pi_V(0,x_0;\mathbb W)_s\right]dh_s.
\end{displaymath}
(Duhamel's principle).
\\
\\
Finally, consider $\mathbb W\in G\Omega_{p,T}(\mathbb R^d)$ and a control $\omega :\bar\Delta_T\rightarrow\mathbb R_+$ satisfying :
\begin{displaymath}
\forall (s,t)\in\bar\Delta_T\textrm{$,$ }
\|\mathbb W\|_{p\textrm{-var};s,t}\leqslant
\omega^{1/p}(s,t).
\end{displaymath}
Then,
\begin{enumerate}
 \item There exists a constant $C_1 > 0$, not depending on $\mathbb W$ and $\omega$, such that for every continous function $h : [0,T]\rightarrow\mathbb R^d$ of finite $q$-variation,
 \begin{displaymath}
 \|D_{h}^{\textrm{FV}}\pi_V(0,x_0;\mathbb W)\|_{\infty;T}
 \leqslant
 C_1\exp\left[C_1(\|h\|_{q\textrm{-var};T}^{p} +
 M_{1,I,\omega})\right].
 \end{displaymath}
 \item There exists a constant $C_2 > 0$, not depending on $\mathbb W$ and $\omega$, such that for every continous function $h : [0,T]\rightarrow\mathbb R^d$ of finite $q$-variation,
 \begin{displaymath}
 \|D_{h}^{\textrm{FV}}\pi_V(0,x_0;\mathbb W)\|_{p\textrm{-var};T}
 \leqslant
 C_2\exp\left[C_2\left[\|h\|_{q\textrm{-var};T}^{p} + \omega(0,T)\right]\right].
 \end{displaymath}
\end{enumerate}

\end{theorem}
%


%
\begin{proof}
Refer to \cite{FV08}, theorems 11.3-6 and Exercice 11.9 for a proof of the first part.
\\
\\
Consider a continuous function $h : [0,T]\rightarrow\mathbb R^d$ of finite $q$-variation, $\mathbb W^h := S_{[p]}(\mathbb W\oplus h)$, $a := (x_0,0,0)$ and
\begin{displaymath}
X_0 :=
\left(1,a,\dots,\frac{a^{\otimes [p]}}{[p]!}\right)\in T^{[p]}\left(\mathbb R^{e + 2}\right).
\end{displaymath}
By \cite{FV08}, Theorem 11.3 :
\begin{displaymath}
D_{h}^{\textrm{FV}}\pi_V(0,x_0;\mathbb W) =
\pi_A(0,0;.)\circ
\mathbb J(.,F)\circ
\pi_G(0,X_0;.)(\mathbb W^h)
\end{displaymath}
where,
\begin{eqnarray*}
 A :\mathbb R^e
 & \longrightarrow &
 \mathcal L(\mathcal L(\mathbb R^e)\times\mathbb R^e ; \mathbb R^e),\\
 F :\mathbb R^e\times\mathbb R^d\times\mathbb R^d
 & \longrightarrow &
 \mathcal L(\mathbb R^e\times\mathbb R^d\times\mathbb R^d;\mathcal L(\mathbb R^e)\times\mathbb R^e)
 \textrm{ and}\\
 G :\mathbb R^e
 & \longrightarrow &
 \mathcal L(\mathbb R^d\times\mathbb R^d;\mathbb R^e\times\mathbb R^d\times\mathbb R^d)
\end{eqnarray*}
are three collections of vector fields, respectively defined by :
\begin{eqnarray*}
 A(a)(L,b) & := & L.a + b,\\
 F(a,a',a'')(b,b',b'') & := & (\langle DV(a),.\rangle b' ; V(a)b'')
 \textrm{ and}\\
 G(a)(b',b'') & := & (V(a)b',b',b'')
\end{eqnarray*}
for every $a,b\in\mathbb R^e$, $a',b',a'',b''\in\mathbb R^d$ and $L\in\mathcal L(\mathbb R^e)$.
\\
\\
By applying successively \cite{FV08}, theorems 10.47 and 10.36, for every $(s,t)\in\bar\Delta_T$,
\begin{eqnarray*}
 \omega_{1}^{1/p}(s,t) & := &
 \left\|\int F\left[\pi_G(0,X_0;\mathbb W^h)\right]d\pi_G(0,X_0;\mathbb W^h)\right\|_{p\textrm{-var};s,t}\\
 & \leqslant &
 \omega_{2}^{1/p}(s,t)
\end{eqnarray*}
with
\begin{displaymath}
\omega_{2}^{1/p}(s,t) :=
\omega_{3}^{1/p}(s,t)\vee
\omega_3(s,t)\vee
\omega_{3}^{p}(s,t)
\end{displaymath}
and, by \cite{FV08}, Proposition 7.52 :
\begin{displaymath}
\omega_3(s,t) :=
\varepsilon_1\left[\|h\|_{q\textrm{-var};s,t}^{p} + \omega(s,t)\right]
\geqslant
\varepsilon_2\|\mathbb W^h\|_{p\textrm{-var};s,t}^{p}
\end{displaymath}
where, $\varepsilon_1,\varepsilon_2\geqslant 1$ are two constants not depending on $\mathbb W$, $\omega$ and $h$.
\\
\\
On one hand, by \cite{FV08}, Exercice 10.55, there exists a constant $C_3 > 0$, not depending on $\mathbb W$, $\omega$ and $h$, such that :
\begin{eqnarray*}
 \left\|D_{h}^{\textrm{FV}}\pi_V(0,x_0;\mathbb W)\right\|_{\infty;T} & \leqslant &
 C_3\exp\left[C_3
 \sup_{
 \begin{tiny}
 \begin{array}{rcl}
  D = \{r_k\} & \in & D_I\\
  \omega_2\left(r_k,r_{k + 1}\right) & \leqslant & 1
 \end{array}
 \end{tiny}
 }
 \sum_{k = 1}^{|D| - 1}
 \omega_2(r_k,r_{k + 1})\right]\\
 & = &
 C_3\exp\left[C_3
 \sup_{
 \begin{tiny}
 \begin{array}{rcl}
  D = \{r_k\} & \in & D_I\\
  \omega_3\left(r_k,r_{k + 1}\right) & \leqslant & 1
 \end{array}
 \end{tiny}
 }
 \sum_{k = 1}^{|D| - 1}
 \omega_3(r_k,r_{k + 1})\right]\\
 & \leqslant &
 C_1\exp\left[C_1\left(
 \|h\|_{q\textrm{-var};T}^{p} + M_{1,I,\omega}\right)\right]
\end{eqnarray*}
with $C_1 := C_3\varepsilon_1$, because
\begin{equation}\label{CTRLegalite}
\omega_2\equiv\omega_3
\textrm{ when }
\omega_2\leqslant 1
\end{equation}
and
\begin{displaymath}
\forall (s,t)\in\bar\Delta_T
\textrm{, }
\omega(s,t)\leqslant\omega_3(s,t).
\end{displaymath}
On the other hand, by \cite{FV08}, Theorem 10.53, there exists a constant $C_4 > 0$, not depending on $\mathbb W$, $\omega$ and $h$, such that for every $(s,t)\in\bar\Delta_T$ satisfying $\omega_2(s,t)\leqslant 1$,
\begin{eqnarray*}
 \|D_{h}^{\textrm{FV}}\pi_V(0,x_0;\mathbb W)_{s,t}\| & \leqslant &
 C_4\left[1 + \|D_{h}^{\textrm{FV}}\pi_V(0,x_0;\mathbb W)_s\|\right]\omega_{2}^{1/p}(s,t)e^{C_4\omega_2(s,t)}\\
 & \leqslant &
 C_4\left[1 + \|D_{h}^{\textrm{FV}}\pi_V(0,x_0;\mathbb W)\|_{\infty;T}\right]\omega_{3}^{1/p}(s,t)e^{C_4\omega_3(0,T)},
\end{eqnarray*}
by (\ref{CTRLegalite}).
\\
\\
Therefore, by super-additivity of the control $\omega_3$, there exists a constant $C_2 > 0$, not depending on $\mathbb W$, $\omega$ and $h$, such that :
\begin{displaymath}
\|D_{h}^{\textrm{FV}}\pi_V(0,x_0;\mathbb W)\|_{p\textrm{-var};T}
\leqslant
C_2\exp\left[C_2\left[\|h\|_{q\textrm{-var;T}}^{p} +\omega(0,T)\right]\right].
\end{displaymath}
\end{proof}
\noindent
At the following corollary, the upper-bound provided at the previous theorem are extended to RDEs having a drift term.
%


%
\begin{corollary}\label{MAJDFVdrift}
Consider $m\in\mathbb N^*$, $p > q\geqslant 1$ such that $1/p + 1/q > 1$, $r\in [1,p[$ such that $1/p + 1/r > 1$, $g : [0,T]\rightarrow\mathbb R^m$ a continuous function of finite $r$-variation, $\mathbb W\in G\Omega_{p,T}(\mathbb R^d)$, $\mathbb W^g := S_{[p]}(\mathbb W\oplus g)$, $V := (V_1,\dots,V_{d + m})$ a collection of $\gamma$-Lipschitz vector fields on $\mathbb R^e$ and $x_0\in\mathbb R^e$. Then, there exists a constant $C > 0$, not depending on $g$ and $\mathbb W$, such that for every continuous function $h : [0,T]\rightarrow\mathbb R^{d + m}$ of finite $q$-variation,
\begin{displaymath}
\|D_{h}^{\textrm{FV}}\pi_V(0,x_0;\mathbb W^g)\|_{\infty;T}
\leqslant
C\exp\left[C\left[\|h\|_{q\textrm{-var};T}^{p} + \|g\|_{r\textrm{-var};T}^{p} + M_{1,I,p}(\mathbb W)\right]\right].
\end{displaymath}
\end{corollary}
%


%
\begin{proof}
Let $h : [0,T]\rightarrow\mathbb R^{d + m}$ be a continuous function of finite $q$-variation. By Corollary \ref{MAJcouplage}, there exists a constant $C_1 > 0$, depending only on $p$ and $r$, such \mbox{that :}
\begin{displaymath}
M_{1,I,p}(\mathbb W^g)
\leqslant
C_1\left[\|g\|_{r\textrm{-var};T}^{p} + M_{1,I,p}(\mathbb W)\right].
\end{displaymath}
Then, by Theorem \ref{DFVsignal}, there exists a constant $C_2 > 0$, not depending on $\mathbb W^g$ and $h$, such that :
\begin{eqnarray*}
 \|D_{h}^{\textrm{FV}}\pi_V(0,x_0;\mathbb W^g)\|_{\infty;T}
 & \leqslant &
 C_2\exp\left[C_2\left[\|h\|_{q\textrm{-var};T}^{p} + M_{1,I,p}(\mathbb W^g)\right]\right]\\
 & \leqslant &
 C\exp\left[C\left[\|h\|_{q\textrm{-var};T}^{p} + \|g\|_{r\textrm{-var};T}^{p} + M_{1,I,p}(\mathbb W)\right]\right]
\end{eqnarray*}
with $C := C_2(1\vee C_1)$.
\end{proof}
\noindent
Finally, let study the second order continuous differentiability of the It\^o map with respect to the driving signal. In order to apply Fernique's theorem at Section 3, tailor-maid upper-bounds are provided for $p\in [1,2[$.
%


%
\begin{lemma}\label{LODEdiff}
Consider $p\in [1,2[$, $x_0\in\mathbb R^e$ and $A := (A_1,\dots, A_d)$ the collection of affine vector fields on $\mathbb R^e$ defined by :
\begin{displaymath}
\forall L\in\mathcal L(\mathbb R^e)
\textrm{$,$ }
\forall a,b\in\mathbb R^e
\textrm{$,$ }
A(a)(L,b) :=
L.a + b.
\end{displaymath}
The map $(L,w)\mapsto\pi_A[0,x_0;(L,w)]$ is continuously differentiable from
\begin{displaymath}
C^{p\textrm{-var}}([0,T];\mathcal L(\mathbb R^e)\times\mathbb R^e)
\textrm{ into }
C^{p\textrm{-var}}([0,T];\mathbb R^e).
\end{displaymath}
Moreover, consider a continuous map $(L,w) : [0,T]\rightarrow\mathcal L(\mathbb R^e)\times\mathbb R^e$ of finite $p$-variation, and $\omega :\bar\Delta_T\rightarrow\mathbb R_+$ a control satisfying :
\begin{displaymath}
\forall (s,t)\in\bar\Delta_T\textrm{$,$ }
\|(L,w)\|_{p\textrm{-var};s,t}^{p}
\leqslant
\omega(s,t).
\end{displaymath}
Then,
\begin{enumerate}
 \item There exists a constant $C_1 > 0$, not depending on $w^L := (L,w)$ and $\omega$, such that for every continuous map $h^H := (H,h) : [0,T]\rightarrow\mathcal L(\mathbb R^e)\times\mathbb R^e$ of finite $p$-variation,
 \begin{displaymath}
 \|\partial_{w^L}\pi_A(0,x_0;w^L).h^H\|_{\infty;T}
 \leqslant
 C_1(1 +\|x_0\|)
 \exp\left[C_1(\|h^H\|_{p\textrm{-var};T}^{p} + M_{1,I,\omega})\right].
 \end{displaymath}
 \item There exists a constant $C_2 > 0$, not depending on $w^L := (L,w)$ and $\omega$, such that for every continuous map $h^H := (H,h) : [0,T]\rightarrow\mathcal L(\mathbb R^e)\times\mathbb R^e$ of finite $q$-variation,
 \begin{displaymath}
 \|\partial_{w^L}\pi_A(0,x_0;w^L).h^H\|_{p\textrm{-var};T}
 \leqslant
 C_2(1 + \|x_0\|)\exp\left[C_2\left[\|h^H\|_{p\textrm{-var};T}^{p} + \omega(0,T)\right]\right].
 \end{displaymath}
\end{enumerate}
\end{lemma}
%


%
\begin{proof}
Since $A$ is a collection of affine vector fields on $\mathbb R^e$, they are locally $\gamma$-Lipschitz on $\mathbb R^e$. Then, by \cite{FV08}, Theorem 10.3, $\pi_A(0,x_0;.)$ is derivable at every points and in every directions of $C^{p\textrm{-var}}([0,T];\mathcal L(\mathbb R^e)\times\mathbb R^e)$.
\\
\\
Consider a continuous map $h^H := (H,h) : [0,T]\rightarrow\mathcal L(\mathbb R^e)\times\mathbb R^e$ of finite $p$-variation.
\\
\\
Firstly,
\begin{eqnarray*}
 D_{h^H}\pi_A(0,x_0;w^L) & = &
 \int_{0}^{.}dL_s.D_{h^H}\pi_A(0,x_0;w^L)_s +
 \int_{0}^{.}dH_s.\pi_A(0,x_0;w^L)_s + h\\
 & = &
 \int_{0}^{.}
 A_1\left[\pi_A(0,x_0;w^L)_s;
 D_{h^H}\pi_A(0,x_0;w^L)_s\right]
 (dw_{s}^{L},dh_{s}^{H})
\end{eqnarray*}
where, $A_1 : \mathbb R^e\times\mathbb R^e\rightarrow\mathcal L(\mathcal L(\mathbb R^e)\times\mathbb R^e\times\mathcal L(\mathbb R^e)\times\mathbb R^e;\mathbb R^e)$ is the collection of affine vector fields defined by :
\begin{displaymath}
A_1(a,b)(c,d) :=
c^1.b +
d^1.a + d^2
\end{displaymath}
for every $a,b\in\mathbb R^e$ and $c,d\in\mathcal L(\mathbb R^e)\times\mathbb R^e$.
\\
\\
By putting $Q(w^L,h^H) := (\pi_A(0,x_0;w^L);D_{h^H}\pi_A(0,x_0;w^L))$ :
\begin{displaymath}
Q(w^L,h^H) =
(x_0,0) +
\int_{0}^{.}
A_2\left[Q_s(w^L,h^H)\right](dw_{s}^{L},dh_{s}^{H})
\end{displaymath}
where, $A_2 : \mathbb R^e\times\mathbb R^e\rightarrow\mathcal L(\mathcal L(\mathbb R^e)\times\mathbb R^e\times\mathcal L(\mathbb R^e)\times\mathbb R^e;\mathbb R^e\times\mathbb R^e)$ is the collection of affine vector fields defined by :
\begin{eqnarray*}
 A_2(a,b)(c,d) & := &
 \begin{bmatrix}
 A(a)c\\
 A_1(a,b)(c,d)
 \end{bmatrix}\\
 & = &
 \begin{pmatrix}
 0\\
 c^1.b
 \end{pmatrix}
 +
 \begin{pmatrix}
 c^1.a\\
 d^1.a
 \end{pmatrix}
 +
 \begin{pmatrix}
 c^2\\
 d^2
 \end{pmatrix}
\end{eqnarray*}
for every $a,b\in\mathbb R^e$ and $c,d\in\mathcal L(\mathbb R^e)\times\mathbb R^e$.
\\
\\
Then, $Q(w^L,h^H) = \pi_{A_2}[0,(x_0,0);(w^L,h^H)]$.
\\
\\
Secondly, there exists a constant $C_3\geqslant 1$, depending only on $p$, such that for every $(s,t)\in\bar\Delta_T$,
\begin{displaymath}
\|(w^L,h^H)\|_{p\textrm{-var};s,t}^{p}
\leqslant
\omega_1(s,t) :=
C_3\left[\|h^H\|_{p\textrm{-var};s,t}^{p} +
\omega(s,t)\right].
\end{displaymath}
By construction, $\omega_1 :\bar\Delta_T\rightarrow\mathbb R_+$ is a control.
\\
\\
By \cite{FV08}, Exercice 10.55, there exists a constant $C_4 > 0$, not depending on $w^L$, $\omega$ and $h^H$, such that :
\begin{eqnarray*}
 \|D_{h^H}\pi_A(0,x_0;w^L)\|_{\infty;T} & \leqslant &
 \|Q(w^L,h^H)\|_{\infty;T}\\
 & \leqslant &
 C_4(1 +\|x_0\|)
 \exp\left[C_4
 \sup_{
 \begin{tiny}
 \begin{array}{rcl}
  D = \{r_k\} & \in & D_I\\
  \omega_1\left(r_k,r_{k + 1}\right) & \leqslant & 1
 \end{array}
 \end{tiny}
 }
 \sum_{k = 1}^{|D| - 1}
 \omega_1(r_k,r_{k + 1})\right]\\
 & \leqslant &
 C_1(1 +\|x_0\|)
 \exp[C_1[
 \|h^H\|_{p\textrm{-var};T}^{p} +\\
 & &
 \sup_{
 \begin{tiny}
 \begin{array}{rcl}
  D = \{r_k\} & \in & D_I\\
  \omega\left(r_k,r_{k + 1}\right) & \leqslant & 1
 \end{array}
 \end{tiny}
 }
 \sum_{k = 1}^{|D| - 1}
 \omega(r_k,r_{k + 1})]]
\end{eqnarray*}
with $C_1 = C_3C_4$, because
\begin{displaymath}
\forall (s,t)\in\bar\Delta_T
\textrm{, }
\omega(s,t)\leqslant
\omega_1(s,t).
\end{displaymath}
Thirdly, by \cite{FV08}, Theorem 10.53, there exists a constant $C_5 > 0$, not depending on $w^L$, $\omega$ and $h^H$, such that for every $(s,t)\in\bar\Delta_T$,
\begin{displaymath}
\|D_{h^H}\pi_A(0,x_0;w^L)_{s,t}\|
\leqslant
C_5(1 + \|x_0\|)\omega_{1}^{1/p}(s,t)e^{C_5\omega_1(0,T)}.
\end{displaymath}
Therefore, by super-additivity of the control $\omega_1$, there exists a constant $C_2 > 0$, not depending on $w^L$, $\omega$ and $h^H$, such that :
\begin{displaymath}
\|D_{h^H}\pi_A(0,x_0;w^L)\|_{p\textrm{-var};T}
\leqslant
C_2(1 + \|x_0\|)\exp\left[C_2\left[\|h^H\|_{p\textrm{-var};T}^{p} + \omega(0,T)\right]\right].
\end{displaymath}
Finally, by the remark following Theorem \ref{ITORDEcontinuite}, since $Q =\pi_{A_2}[0,(x_0,0);.]$ and $A_2$ is a collection of affine vector fields on $\mathbb R^e$, $Q$ and (then) $(w^L,h^H)\mapsto D_{h^H}\pi_A(0,x_0;w^L)$ are uniformly continuous on bounded set of
\begin{displaymath}
C^{p\textrm{-var}}([0,T];\mathcal L(\mathbb R^e)\times\mathbb R^e)\times
C^{p\textrm{-var}}([0,T];\mathcal L(\mathbb R^e)\times\mathbb R^e).
\end{displaymath}
In conclusion, by \cite{FV08}, Proposition B.5, $\pi_A(0,x_0;.)$ is continuously differentiable as stated.
\end{proof}
%


%
\begin{proposition}\label{DVITO2diff}
Consider $2 > p\geqslant q\geqslant 1$, $x_0\in\mathbb R^e$ and $V,\tilde V\in\textrm{Lip}^{\gamma + 1}(\mathbb R^e;\mathbb R^d)$. The maps
\begin{eqnarray*}
 D^1 & : &
 \left\{
 \begin{array}{rcl}
 C^{p\textrm{-var}}([0,T];\mathbb R^d) & \longrightarrow & C^{p\textrm{-var}}([0,T];\mathbb R^{e\times e})\\
 w & \longmapsto & D^1(w) = J_{.\leftarrow 0}^{w}
 \end{array}\right.,
 \\
 D^2 & : &
 \left\{
 \begin{array}{rcl}
 C^{p\textrm{-var}}([0,T];\mathbb R^d) & \longrightarrow & C^{p\textrm{-var}}([0,T];\mathbb R^{e\times e})\\
 w & \longmapsto & D^2(w) = J_{0\leftarrow .}^{w}
 \end{array}\right.
 \textrm{ and}\\
 D^3 & : &
 \left\{
 \begin{array}{rcl}
 C^{p\textrm{-var}}([0,T];\mathbb R^d) & \longrightarrow & C^{p\textrm{-var}}([0,T];\mathbb R^e)\\
 w & \longmapsto & D^3(w) = \partial_V\pi_V(0,x_0;w).\tilde V
 \end{array}\right.
\end{eqnarray*}
are continuously differentiable.
\\
\\
Moreover, let $w : [0,T]\rightarrow\mathbb R^d$ be a continuous map of finite $p$-variation. there exists three strictly positive constants $C_1$, $C_2$ and $C_3$, not depending on $w$, such that for every continuous function $h : [0,T]\rightarrow\mathbb R^d$ of finite $q$-variation,
\begin{eqnarray*}
 \|\partial_wD^1(w).h\|_{p\textrm{-var};T} & \leqslant & C_1\exp\left[C_1(\|h\|_{q\textrm{-var};T}^{p} +\|w\|_{p\textrm{-var};T}^{p})\right],\\
 \|\partial_wD^2(w).h\|_{p\textrm{-var};T} & \leqslant & C_2\exp\left[C_2(\|h\|_{q\textrm{-var};T}^{p} +\|w\|_{p\textrm{-var};T}^{p})\right]
 \textrm{ and}\\
 \|\partial_wD^3(w).h\|_{p\textrm{-var};T} & \leqslant & C_3\exp\left[C_3(\|h\|_{q\textrm{-var};T}^{p} +\|w\|_{p\textrm{-var};T}^{p})\right].
\end{eqnarray*}
\end{proposition}
%


%
\begin{proof}
Proofs are similar for $D^1$, $D^2$ and $D^3$. It is detailed for $D^3$.
\\
\\
At the proof of Theorem \ref{CIVFdiff}, with the same notations, it has been established that :
\begin{displaymath}
D^3(w) =
(I\circ J\circ K)(w)
\end{displaymath}
where,
\begin{displaymath}
I := \pi_A(0,0;.)
\textrm{, }
J := \mathbb J(.;F_{V,\tilde V})
\textrm{ and }
K := \pi_{F_V}\left[0,(x_0,0);.\right].
\end{displaymath}
When $p\in [1,2[$, since $F_V$ and $F_{V,\tilde V}$ are respectively two collections of $\gamma$-Lipschitz vector fields on $\mathbb R^e$ and $\mathbb R^e\times\mathbb R^d$ by construction, by Theorem \ref{DFVsignal} :
\begin{eqnarray*}
 J : C^{p\textrm{-var}}([0,T];\mathbb R^e\times\mathbb R^d)
 & \longrightarrow &
 C^{p\textrm{-var}}([0,T];\mathcal L(\mathbb R^e)\times\mathbb R^e)
 \textrm{ and}\\
 K : C^{p\textrm{-var}}([0,T];\mathbb R^d)
 & \longrightarrow &
 C^{p\textrm{-var}}([0,T];\mathbb R^e\times\mathbb R^d)
\end{eqnarray*}
are continuously differentiable.
\\
\\
Moreover, by Proposition \ref{LODEdiff}, $I$ is a continuously differentiable map from
\begin{displaymath}
C^{p\textrm{-var}}([0,T];\mathcal L(\mathbb R^e)\times\mathbb R^e)
\textrm{ into }
C^{p\textrm{-var}}([0,T];\mathbb R^e).
\end{displaymath}
Then, by composition, $D^3$ is a continuously differentiable map from
\begin{displaymath}
C^{p\textrm{-var}}([0,T];\mathbb R^d)
\textrm{ into }
C^{p\textrm{-var}}([0,T];\mathbb R^e).
\end{displaymath}
\noindent
Let $h : [0,T]\rightarrow\mathbb R^d$ be a continuous function of finite $q$-variation :
\begin{displaymath}
\partial_wD^3(w).h =
\langle DI[(J\circ K)(w)];
\langle DJ[K(w)],\partial_wK(w).h\rangle\rangle.
\end{displaymath}
Then,
\begin{eqnarray}
 \label{DVITO2diffmaj}
 \|\partial_wD^3(w).h\|_{p\textrm{-var};T} & \leqslant &
 \|DI[(J\circ K)(w)]\|_{\mathcal L;p}\times\\
 & &
 \|DJ[K(w)]\|_{\mathcal L;p}\times
 \nonumber\\
 & &
 \|\partial_w K(w).h\|_{p\textrm{-var};T}
 \nonumber
\end{eqnarray}
where, $\|.\|_{\mathcal L;p}$ is the usual operator norm on
\begin{displaymath}
\mathcal L(
C^{p\textrm{-var}}([0,T];\mathcal L(\mathbb R^e)\times\mathbb R^e) ;
C^{p\textrm{-var}}([0,T];\mathbb R^e))
\end{displaymath}
or on
\begin{displaymath}
\mathcal L(
C^{p\textrm{-var}}([0,T];\mathbb R^e\times\mathbb R^d) ;
C^{p\textrm{-var}}([0,T];\mathcal L(\mathbb R^e)\times\mathbb R^e)).
\end{displaymath}
Let's find upper-bounds for each term of the product of the right-hand side of inequality (\ref{DVITO2diffmaj}) :
\begin{enumerate}
 \item Since $K(w) =\pi_{F_V}(0,x_0;w)$, by the second point of Theorem \ref{DFVsignal} for the control $\omega_{w,p}$, there exists a constant $C_4 > 0$, not depending on $w$ and $h$, such that :
 \begin{equation}\label{DVITO2diffmaj1}
 \|\partial_wK(w).h\|_{p\textrm{-var};T}
 \leqslant
 C_4\exp\left[C_4(\|h\|_{q\textrm{-var};T}^{p} +
 \|w\|_{p\textrm{-var};T}^{p})\right].
 \end{equation}
 \item Let $g : [0,T]\rightarrow\mathbb R^e\times\mathbb R^d$ be a continuous function of finite $p$-variation such that $\|g\|_{p\textrm{-var};T}\leqslant 1$.
 \\
 \\
 On one hand, by construction of the rough integral :
 \begin{displaymath}
 (J\circ K)(w) =
 \pi_{\mathcal L(\mathbb R^e)\times\mathbb R^e}[
 \pi_{G_{V,\tilde V}}[0,(x_0,0);K(w)]]
 \end{displaymath}
 where,
 \begin{displaymath}
 G_{V,\tilde V} : (\mathbb R^e\times\mathbb R^d)^2
 \longrightarrow
 \mathcal L(\mathbb R^e\times\mathbb R^d;
 \mathbb R^d\times\mathbb R^e\times
 \mathcal L(\mathbb R^e)\times\mathbb R^e)
 \end{displaymath}
 is the collection of $\gamma$-Lipschitz vector fields (as $F_{V,\tilde V}$) defined by :
 \begin{displaymath}
 \forall a,b,c\in
 \mathbb R^e\times\mathbb R^d
 \textrm{, }
 G_{V,\tilde V}(a,b)c =
 (c,F_{V,\tilde V}(a).c).
 \end{displaymath}
 On the other hand, by \cite{FV08}, Theorem 10.36, for every $(s,t)\in\bar\Delta_T$,
 \begin{displaymath}
 \|K(w)\|_{p\textrm{-var};s,t}
 \leqslant
 \omega_{1}^{1/p}(s,t) :=
 \omega_{2}^{1/p}(s,t)\vee\omega_2(s,t)
 \end{displaymath}
 with
 \begin{displaymath}
 \omega_2(s,t) :=
 \varepsilon_1\|w\|_{p\textrm{-var};s,t}^{p}
 \end{displaymath}
 where, $\varepsilon_1\geqslant 1$ is a constant not depending on $w$. However, for every $\varepsilon\geqslant 1$,
 \begin{equation}\label{DFVsignalpart}
 \omega_3(.;\varepsilon)
 \equiv
 \varepsilon(
 \|g\|_{p\textrm{-var}}^{p} +
 \omega_2)
 \textrm{ when }
 \omega_3(.;\varepsilon)\leqslant 1
 \end{equation}
 where, for every $(s,t)\in\bar\Delta_T$,
 \begin{displaymath}
 \omega_{3}^{1/p}(s,t;\varepsilon) :=
 \omega_{4}^{1/p}(s,t;\varepsilon)\vee
 \omega_4(s,t;\varepsilon)\vee
 \omega_{4}^{p}(s,t;\varepsilon)
 \end{displaymath}
 and
 \begin{displaymath}
 \omega_4(s,t;\varepsilon) :=
 \varepsilon\left[
 \|g\|_{p\textrm{-var};s,t}^{p} +
 \omega_1(s,t)\right].
 \end{displaymath}
 Then, as at the proof of point 2 of Theorem \ref{DFVsignal}, by using (\ref{DFVsignalpart}) ; there exists a constant $C_ 5 > 0$, not depending on $w$ and $g$, such that :
 \begin{displaymath}
 \|DJ[K(w)].g\|_{p\textrm{-var};T}
 \leqslant
 C_5\exp[C_5(\|g\|_{p\textrm{-var};T}^{p} +
 \|w\|_{p\textrm{-var};T}^{p})].
 \end{displaymath}
 Therefore,
 \begin{equation}\label{DVITO2diffmaj2}
 \|DJ[K(w)]\|_{\mathcal L;p}
 \leqslant
 C_5\exp[C_5(1 +
 \|w\|_{p\textrm{-var};T}^{p})].
 \end{equation}
 \item Let $g^H := (H,g) : [0,T]\rightarrow\mathcal L(\mathbb R^e)\times\mathbb R^e$ be the continuous map of finite $p$-variation such that $\|g\|_{p\textrm{-var};T}\leqslant 1$.
 \\
 \\
 On one hand, it has been established at Lemma \ref{LODEdiff} that there exists a collection $A_1$ of affine vector fields on $\mathbb R^{e_1}\times\mathbb R^{e_2}$ ($e_1 = e_2 = e$) such that :
 \begin{displaymath}
 DI[(J\circ K)(w)].g^H =
 \pi_{\mathbb R^{e_2}}[
 \pi_{A_1}[0,(0,0);((J\circ K)(w),g^H)].
 \end{displaymath}
 On the other hand, by applying successively \cite{FV08}, theorems 10.36 and 10.47, for every $(s,t)\in\bar\Delta_T$,
 \begin{displaymath}
 \|(J\circ K)(w)\|_{p\textrm{-var};s,t}
 \leqslant
 \omega_{5}^{1/p}(s,t) :=
 \omega_{6}^{1/p}(s,t)\vee
 \omega_6(s,t)\vee
 \omega_{6}^{p}(s,t)
 \end{displaymath}
 with
 \begin{displaymath}
 \omega_6(s,t) :=
 \varepsilon_2\|w\|_{p\textrm{-var};s,t}^{p}
 \end{displaymath}
 where, $\varepsilon_2\geqslant 1$ is a constant not depending on $w$.
 \\
 \\
 Then, as at the proof of point 2 of Lemma \ref{LODEdiff} ; there exists a constant $C_6 > 0$, not depending on $w$ and $g^H$, such that :
 \begin{displaymath}
 \|DI[(J\circ K)(w)].g^H\|_{p\textrm{-var};T}
 \leqslant
 C_6\exp[C_6(\|g^H\|_{p\textrm{-var};T}^{p} +
 \|w\|_{p\textrm{-var};T}^{p})].
 \end{displaymath}
 Therefore,
 \begin{equation}\label{DVITO2diffmaj3}
 \|DI[(J\circ K)(w)]\|_{\mathcal L;p}
 \leqslant
 C_6\exp[C_6(1 +
 \|w\|_{p\textrm{-var};T}^{p})].
 \end{equation}
\end{enumerate}
In conclusion, inequalities (\ref{DVITO2diffmaj1}), (\ref{DVITO2diffmaj2}), (\ref{DVITO2diffmaj3}) and (\ref{DVITO2diffmaj}) imply together that there exists a constant $C_3 > 0$, not depending on $w$ and $h$, such that :
\begin{displaymath}
\|\partial_wD^3(w).h\|_{p\textrm{-var};T}
\leqslant
C_3\exp\left[C_3(\|h\|_{q\textrm{-var};T}^{p} +\|w\|_{p\textrm{-var};T}^{p})\right].
\end{displaymath}
\end{proof}
%


%
\subsection{Application to Gaussian stochastic analysis}
Consider a $d$-dimensional stochastic process $W$ and the probability space $(\Omega,\mathcal A,\mathbb P)$, where $\Omega$ is the vector space of continuous functions from $[0,T]$ into $\mathbb R^d$, $\mathcal A$ is the $\sigma$-algebra generated by cylinder sets of $\Omega$, and $\mathbb P$ is the probability measure induced by the process $W$ on $(\Omega,\mathcal A)$.
\\
\\
In order to prove Corollary \ref{ITODPint} that is crucial at Section 3, let's first state existing results on Gaussian rough paths established by P. Friz and N. Victoir in \cite{FV10}, and by T. Cass, C. Litterer and T. Lyons in \cite{CLL11}.
\\
\\
Consider the two following technical assumptions on the stochastic process $W$.
%


%
\begin{assumption}\label{CSGRPgauss}
$W$ is a $d$-dimensionnel centered Gaussian process with continuous paths. Moreover, its covariance function $c_W$ is of finite 2D $\rho$-variation with $\rho\in [1,2[$ (cf. \cite{FV08}, Definition 5.50).
\end{assumption}
%


%
\begin{assumption}\label{CMhyp}
There exists $p > q\geqslant 1$ such that :
\begin{displaymath}
\frac{1}{p} +\frac{1}{q} > 1
\textrm{ and }
H^1\hookrightarrow
C^{q\textrm{-var}}([0,T];\mathbb R^d).
\end{displaymath}
\end{assumption}
\noindent
\textbf{Example.} By \cite{FV08}, Proposition 15.5, Proposition 15.7 and Exercice 20.2, a fractional Brownian motion of Hurst parameter $H\in ]1/4,1/2]$ satisfies assumptions \ref{CSGRPgauss} and \ref{CMhyp}.
%


%
\begin{theorem}\label{GRPgauss}
Consider a stochastic process $W$ satisfying Assumption \ref{CSGRPgauss}, and $p > 2\rho$. For almost every $\omega\in\Omega$, there exists a geometric $p$-rough path $\mathbb W(\omega)$ over $W(\omega)$ satisfying :
\begin{enumerate}
 \item There exists a deterministic constant $C > 0$, only depending on $\rho$, $p$ and $\|c_W\|_{\rho\textrm{-var};[0,T]^2}$, such that :
 \begin{displaymath}
 \mathbb E\left(e^{C\|\mathbb W\|_{p\textrm{-var};T}^{2}}\right) < \infty.
 \end{displaymath}
 (generalized Fernique's theorem).
 \item Let $(W^n,n\in\mathbb N)$ be a sequence of linear approximations, or of mollifier approximations, of the process $W$. Then, $\mathbb W$ is the limit in $p$-variation, in $L^r(\Omega)$ for every $r\geqslant 1$, of the sequence $(S_3(W^n),n\in\mathbb N)$ (universality).
\end{enumerate}
$\mathbb W$ is the enhanced Gaussian process associated to $W$.
\end{theorem}
\noindent
Refer to \cite{FV08}, Theorem 15.33 for a proof.
%


%
\begin{proposition}\label{TRANSgauss}
Consider a stochastic process $W$ satisfying assumptions \ref{CSGRPgauss} and \ref{CMhyp}, $\mathbb W$ the enhanced Gaussian process associated to $W$, and the Cameron-Martin's space $H^1\subset\Omega$ of the process $W$. Then,
\begin{displaymath}
\forall\omega\in\Omega\textrm{$,$ }
\forall h\in H^1\textrm{$,$ }
\mathbb W(\omega + h) =
T_h\mathbb W(\omega).
\end{displaymath}
\end{proposition}
\noindent
Refer to \cite{FV08}, Lemma 15.58 for a proof.
%


%
\begin{proposition}\label{CLLrelation}
For every geometric $p$-rough path $Y$ and every $\varepsilon > 0$,
\begin{displaymath}
M_{\varepsilon,I,p}(Y)
\leqslant
\varepsilon\left[
2N_{\varepsilon,I,p}(Y) + 1
\right].
\end{displaymath}
\end{proposition}
\noindent
Refer to \cite{CLL11}, Proposition 4.6 for a proof.
%


%
\begin{theorem}\label{CLLprincipal}
Consider a stochastic process $W$ satisfying assumptions \ref{CSGRPgauss} and \ref{CMhyp}, and $\mathbb W$ the enhanced Gaussian process associated to $W$. Then,
\begin{displaymath}
\forall C,\varepsilon,r > 0
\textrm{$,$ }
Ce^{CN_{\varepsilon,I,p}(\mathbb W)}\in L^r(\Omega).
\end{displaymath}
\end{theorem}
\noindent
Refer to \cite{CLL11}, Theorem 6.4 and Remark 6.5 for a proof.
%


%
\begin{corollary}\label{ITODPint}
Consider $x_0\in\mathbb R^e$, $V := (V_1,\dots,V_{d + 1})$ and $\tilde V := (\tilde V_1,\dots,\tilde V_{d + 1})$ two collections of $\gamma$-Lipschitz vector fields on $\mathbb R^e$, a stochastic process $W$ satisfying assumptions \ref{CSGRPgauss} and \ref{CMhyp}, $\mathbb W$ the enhanced Gaussian process associated to $W$, $\mathbb W^g := S_{[p]}(\mathbb W\oplus g)$ with $g := \textrm{Id}_{[0,T]}$, and a continuous function $h : [0,T]\rightarrow\mathbb R^{d + 1}$ of finite $q$-variation. Then,
\begin{displaymath}
\|J_{.\leftarrow 0}^{\mathbb W^g}\|_{\infty;T}
\textrm{$,$ }
\|\partial_V\pi_V(0,x_0;\mathbb W^g).\tilde V\|_{\infty;T}
\textrm{ and }
\|D_{h}^{\textrm{FV}}\pi_V(0,x_0;\mathbb W^g)\|_{\infty;T}
\end{displaymath}
belong to $L^r(\Omega)$ for every $r > 0$.
\end{corollary}
%


%
\begin{proof}
It is a straightforward consequence of corollaries \ref{MAJsdrift} and \ref{MAJDFVdrift}, of Proposition \ref{CLLrelation} (deterministic results), and of Theorem \ref{CLLprincipal} (probabilistic result).
\end{proof}
%


%
\section{Sensitivity analysis for Gaussian rough differential equations}
\noindent
This section solves the problematic stated in introduction of the paper by using the deterministic results on RDEs of subsections 2.1 and 2.2, the probabilistic results on Gaussian RDEs of Subsection 2.3 and Malliavin calculus arguments.
\\
\\
Assume that $W$, $\mu$ and $\sigma$ defined in introduction satisfy the following condition.
%


%
\begin{assumption}\label{HYPsynth}
The process $W$ satisfies assumptions \ref{CSGRPgauss} and \ref{CMhyp}, and
\begin{displaymath}
C_{0}^{1}\left([0,T];\mathbb R^d\right)
\subset H^1.
\end{displaymath}
Moreover, there exists a constant $C > 0$ such that :
\begin{displaymath}
\forall h\in C_{0}^{1}\left([0,T];\mathbb R^d\right)
\textrm{$,$ }
\|h\|_{H^1}
\leqslant
C\|\dot h\|_{\infty;T}.
\end{displaymath}
Functions $\mu$ and $\sigma$ satisfy Assumption \ref{REGcoef} and, for every $a\in\mathbb R^d$, $\sigma(a)$ is an invertible matrix. Moreover, the function $\sigma^{-1} : \mathbb R^d\rightarrow\mathcal M_d(\mathbb R)$ is bounded.
\end{assumption}
\noindent
\textbf{Example.} A fractional Brownian motion $B^H$ of Hurst parameter $H\in ]1/4,1[$ satisfies Assumption \ref{HYPsynth}. Indeed, it as been stated at Subsection 2.3 that $B^H$ satisfies assumptions \ref{CSGRPgauss} and \ref{CMhyp}. Moreover, by the first point of L. Decreusefond and S. Ustunel \cite{DU98}, \mbox{Theorem 3.3 :}
\begin{displaymath}
C_{0}^{1}\left([0,T];\mathbb R^d\right)\subset H^1.
\end{displaymath}
Consider $h\in C_{0}^{1}([0,T];\mathbb R^d)$. By the second point of \cite{DU98}, Theorem 3.3 :
\begin{eqnarray*}
 \|h\|_{H^1} & = &
 \|J_H(\dot h)\|_{H^1} =
 \|\dot h\|_{L^2([0,T])}\\
 & \leqslant &
 T^{1/2}\|\dot h\|_{\infty;T}.
\end{eqnarray*}
Assume also that the function $F :\mathbb R^d\rightarrow\mathbb R$ satisfies one of the two following assumptions.
%


%
\begin{assumption}\label{Freg1}
The function $F$ is continuously differentiable from $\mathbb R^d$ into $\mathbb R$. Moreover, there exists two constants $C > 0$ and $N\in\mathbb N^*$ such that, for every $a\in\mathbb R^d$,
\begin{displaymath}
|F(a)|\leqslant
C(1 +\|a\|)^N
\textrm{ and }
\|DF(a)\|_{\mathcal L}
\leqslant
C(1 + \|a\|)^N.
\end{displaymath}
\end{assumption}
%


%
\begin{assumption}\label{Freg2}
There exists two constants $C > 0$ and $N\in\mathbb N^*$ such that, for every $a\in\mathbb R^d$,
\begin{displaymath}
|F(a)|\leqslant
C(1 +\|a\|)^N.
\end{displaymath}
\end{assumption}
\noindent
The following results are solving, at least partially, the problematic stated in introduction of the paper.
\\
\\
\textbf{Notations :}
\begin{itemize}
 \item Under Assumption \ref{HYPsynth}, the enhanced Gaussian process associated to $W$ is denoted by $\mathbb W$, $\mathbb W^g := S_{[p]}(\mathbb W\oplus g)$ with $g := \textrm{Id}_{[0,T]}$, and $V := (V_1,\dots,V_{d + 1})$ is the collection of vector fields defined by :
 \begin{displaymath}
 V(a)(b,c) :=
 \mu(a)c +
 \sigma(a)b
 \end{displaymath}
 for every $a,b\in\mathbb R^d$ and $c\in\mathbb R$.
 \item Let $\mathbb S_p\subset\textrm{Lip}^{\gamma}(\mathbb R^d;\mathbb R^d)$ be the space of functions from $\mathbb R^d$ into $\mathcal M_d(\mathbb R)$, $[p] + 1$ times differentiable, bounded and of bounded derivatives.
 \item For every $x\in\mathbb R^d$, $\mathbb E[F(X_T)]$ is denoted by $f_T(x,\sigma)$.
\end{itemize}
%


%
\begin{lemma}\label{CMiso}
Let $I = (I^1,\dots,I^d)$ be the map from $H$ into $H^1$ such that :
\begin{displaymath}
I^i(h) :=
\mathbb E\left[\mathbf W^i(h^i)W^i\right]\in H^1
\end{displaymath}
for every $h\in H = H_1\oplus\dots\oplus H_d$ and $i = 1,\dots,d$. Then, $I$ is an isometry from $H$ into $H^1$.
\end{lemma}
%


%
\begin{proof}
On one hand, the linearity of $I$ as a map from $H$ into $H^1$ is a straightforward consequence of the linearity of $\mathbf W$ as a map from $H$ into $L^2(\Omega)$.
\\
\\
On the other hand, by construction of $\mathbf W$ and of the scalar products on $H$ and \mbox{$H^1$ :}
\begin{eqnarray*}
 \langle I(h),I(g)\rangle_{H^1} & = &
 \sum_{i = 1}^{d}\langle\mathbb E\left[\mathbf W^i(h^i)W^i\right] ;
 \mathbb E\left[\mathbf W^i(g^i)W^i\right]\rangle_{H_{i}^{1}}\\
 & = &
 \sum_{i = 1}^{d}\mathbb E\left[
 \mathbf W^i(h^i)\mathbf W^i(g^i)
 \right] =
 \langle h,g\rangle_H
\end{eqnarray*}
for every functions $h,g\in H$.
\end{proof}
\noindent
The following corollary extends \cite{FV08}, Proposition 20.5 to Gaussian RDEs having a drift term.
%


%
\begin{lemma}\label{H1DIFFrde}
For every $x_0\in\mathbb R^d$ and almost every $\omega\in\Omega$, the map $h\mapsto\pi_V[0,x_0;\mathbb W^g(\omega + h)]$ is continuously differentiable from
\begin{displaymath}
H^1
\textrm{ into }
C^{p\textrm{-var}}\left([0,T];\mathbb R^d\right).
\end{displaymath}
In particular, for every $t\in [0,T]$, $\pi_V(0,x_0;\mathbb W^g)_t\in\mathbb D_{\textrm{loc}}^{1,2}$ and for every $h\in H^1$,
\begin{eqnarray*}
 \langle\mathbf D\pi_V(0,x_0;\mathbb W^g)_t, I^{-1}(h)\rangle_H & = &
 D_{(h,0)}^{\textrm{FV}}
 \pi_V(0,x_0;\mathbb W^g)_t\\
 & = &
 \int_{0}^{t}
 J_{t\leftarrow s}^{\mathbb W^g}\sigma\left[\pi_V(0,x_0;\mathbb W^g)_s\right]dh_s.
\end{eqnarray*}
\end{lemma}
%


%
\begin{proof}
By Proposition \ref{TRANSgauss}, for almost every $\omega\in\Omega$ and every $h\in H^1$,
\begin{eqnarray*}
 \mathbb W^g(\omega + h) & = &
 S_{[p]}[\mathbb W(\omega + h)\oplus g]\\
 & = &
 S_{[p]}[T_h\mathbb W(\omega)\oplus g]\\
 & = &
 T_{(h,0)}S_{[p]}[\mathbb W(\omega)\oplus g]\\
 & = &
 T_{(h,0)}\mathbb W^g(\omega).
\end{eqnarray*}
Then, almost surely :
\begin{equation}\label{EGH1diff}
\pi_V[0,x_0;\mathbb W^g(. + h)] =
\pi_V[0,x_0;T_{(h,0)}\mathbb W^g].
\end{equation}
Moreover, Assumpution \ref{CMhyp} and Corollary \ref{MAJDFVdrift} imply that $h\mapsto\pi_V[0,x_0;T_{(h,0)}\mathbb W^g]$ is continuously differentiable from
\begin{displaymath}
H^1\subset C^{q\textrm{-var}}\left([0,T];\mathbb R^d\right)\textrm{ into }
C^{p\textrm{-var}}\left([0,T];\mathbb R^d\right).
\end{displaymath}
Then, by equality (\ref{EGH1diff}), the map $h\mapsto\pi_V[0,x_0;\mathbb W^g(. + h)]$ is also continuously differentiable from
\begin{displaymath}
H^1\textrm{ into }
C^{p\textrm{-var}}\left([0,T];\mathbb R^d\right),
\end{displaymath}
and for almost every $\omega\in\Omega$ and every $h\in H^1$,
\begin{displaymath}
D_{(h,0)}^{\textrm{FV}}\pi_V\left[0,x_0;\mathbb W^g(\omega)\right] =
D_hF^{\omega}(0)
\end{displaymath}
with $F^{\omega} :=\pi_V[0,x_0;\mathbb W^g(\omega + .)]$.
\\
\\
Moreover, by Duhamel's principle (cf. Theorem \ref{DFVsignal}), for every $t\in [0,T]$ and every $h\in H^1$,
\begin{eqnarray*}
 D_{(h,0)}^{\textrm{FV}}
 \pi_V(0,x_0;\mathbb W^g)_t
 & = &
 \int_{0}^{t}
 J_{t\leftarrow s}^{\mathbb W^g}V\left[\pi_V(0,x_0;\mathbb W^g)_s\right](dh_s,0)\\
 & = &
 \int_{0}^{t}
 J_{t\leftarrow s}^{\mathbb W^g}\sigma\left[\pi_V(0,x_0;\mathbb W^g)_s\right]dh_s.
\end{eqnarray*}
In conclusion, by \cite{NUALART06}, Proposition 4.1.3 and Lemma 4.1.2, $\pi_V(0,x_0;\mathbb W^g)_t$ is continuously $H^1$-differentiable and then locally derivable in the sense of Malliavin, \mbox{with :}
\begin{displaymath}
\langle\mathbf D\pi_V(0,x_0;\mathbb W^g)_t, I^{-1}(h)\rangle_H =
D_hF^{\omega}(0).
\end{displaymath}
\end{proof}
%


%
\begin{theorem}\label{CALCsensibilites}
Under assumptions \ref{HYPsynth} and \ref{Freg1} :
\begin{enumerate}
 \item The function $f_T(.,\sigma)$ is differentiable from $\mathbb R^d$ into $\mathbb R$ and, for every $x,v\in\mathbb R^d$,
 \begin{displaymath}
 \partial_xf_T(x,\sigma).v =
 \mathbb E\left[\langle
 \mathbf D(F\circ X_{T}^{x}),I^{-1}(h^{x,v})\rangle_H\right]
 \end{displaymath}
where
 \begin{displaymath}
 h^{x,v} :=
 \frac{1}{T}
 \int_{0}^{.}
 \sigma^{-1}\left(X_{s}^{x}\right)J_{s\leftarrow 0}^{\mathbb W^g}vds.
 \end{displaymath}
 \item For every $x\in\mathbb R^d$, the function $f_T(x,.)$ is differentiable from $\mathbb S_p$ into $\mathbb R$ and, for every $\sigma,\tilde\sigma\in\mathbb S_p$ satisfying Assumption \ref{HYPsynth},
 \begin{displaymath}
 \partial_{\sigma}f_T(x,\sigma).\tilde\sigma =
 \mathbb E\left[\langle
 \mathbf D(F\circ X_{T}^{\sigma}),I^{-1}(h^{\sigma,\tilde\sigma})\rangle_H\right]
 \end{displaymath}
 where
 \begin{displaymath}
 h^{\sigma,\tilde\sigma} :=
 \frac{1}{T}
 \int_{0}^{.}
 \sigma^{-1}\left(X_{s}^{\sigma}\right)J_{s\leftarrow T}^{\mathbb W^g}\partial_{\sigma}X_{T}^{\sigma}.\tilde\sigma ds.
 \end{displaymath}
\end{enumerate}
\end{theorem}
%


%
\begin{proof}
Proofs of points 1 and 2 are similar :
\begin{enumerate}
 \item On one hand, for every $\varepsilon\in ]0,1]$, $\eta > 0$ and $x,v\in\mathbb R^d$, by Taylor's formula, and the first point of Corollary \ref{MAJsdrift} ; there exists a constant $C_1 > 0$ depending only on $p$, $\gamma$, $\eta$ and $\|V\|_{\textrm{lip}^{\gamma}}$ such that :
 \begin{eqnarray*}
  \frac{\left|F(X_{T}^{x +\varepsilon v}) - F(X_{T}^{x})\right|}{\varepsilon}
  & = &
  \left|\int_{0}^{1}
  \langle
  DF(X_{T}^{x +\theta\varepsilon v}),
  DX_{T}^{x +\theta\varepsilon v}.v
  \rangle d\theta
  \right|\\
  & \leqslant &
  C_1\|v\|e^{C_1M_{\eta,I,p}(\mathbb W)}
  \int_{0}^{1}
  \|DF(X_{T}^{x +\theta\varepsilon v})\|_{\mathcal L}d\theta.
 \end{eqnarray*}
 Moreover, since $F$ satisfies Assumption \ref{Freg1}, there exists two constants $C_2 > 0$ and $N\in\mathbb N^*$, depending only on $F$, such that for every $\theta\in [0,1]$,
 \begin{displaymath}
 \left\|DF\left(X_{T}^{x +\theta\varepsilon v}\right)\right\|_{\mathcal L}
 \leqslant
 C_2\left(1 +\left\|X_{T}^{x +\theta\varepsilon v}\right\|\right)^N.
 \end{displaymath}
 Then, by the triangle inequality together with \cite{FV08}, Theorem 10.36, there exists a constant $C_3 > 0$, not depending on $x$, $v$, $\theta$, $\varepsilon$, $V$ and $\mathbb W$, such that :
 \begin{eqnarray*}
 \|DF(X_{T}^{x +\theta\varepsilon v})\|_{\mathcal L}
 & \leqslant &
 C_2[1 + \|x\| + \|v\| +\\
 & &
 C_3[\|V\|_{\textrm{lip}^{\gamma - 1}}\|\mathbb W^g\|_{p\textrm{-var};T}\vee
 \|V\|_{\textrm{lip}^{\gamma - 1}}^{p}\|\mathbb W^g\|_{p\textrm{-var};T}^{p}]]^N.
 \end{eqnarray*}
 Since $W$ satisfies assumptions \ref{CSGRPgauss} and \ref{CMhyp}, by Corollary \ref{ITODPint}, the generalized Fernique's theorem (cf. Theorem \ref{GRPgauss}) and Cauchy-Schwarz's inequality :
 \begin{displaymath}
 \varepsilon\in ]0,1]
 \longmapsto
 \frac{\left|F(X_{T}^{x +\varepsilon v}) - F(X_{T}^{x})\right|}{\varepsilon}
 \end{displaymath}
 is bounded by an integrable random variable not depending on $\varepsilon$. Therefore, by Lebesgue's theorem, $f_T(.,\sigma)$ is differentiable on $\mathbb R^d$ and
 \begin{equation}\label{FDIFdelta}
 \forall x,v\in\mathbb R^d
 \textrm{, }
 \partial_x f_T(x,\sigma).v =
 \mathbb E\left[\langle DF(X_{T}^{x}),\partial_x X_{T}^{x}.v\rangle\right].
 \end{equation}
 On the other hand, consider $x,v\in\mathbb R^d$. By construction, paths of the process $h^{x,v}$ are continuously differentiable from $[0,T]$ into $\mathbb R^d$ and $h_{0}^{x,v} = 0$. Then, since $W$ satisfies Assumption \ref{HYPsynth}, $h^{x,v}$ is a $H^1$-valued random variable. By Duhamel's principle (cf. Theorem \ref{DFVsignal}) :
 \begin{eqnarray*}
  D_{(h^{x,v},0)}^{\textrm{FV}}X_{T}^{x} & = &
  \int_{0}^{T}
  J_{T\leftarrow s}^{\mathbb W^g}\sigma\left(X_{s}^{x}\right)dh_{s}^{x,v}\\
  & = &
  \partial_x X_{T}^{x}.v.
 \end{eqnarray*}
 Therefore, by equality (\ref{FDIFdelta}), Lemma \ref{H1DIFFrde} and \cite{NUALART06}, Proposition 1.2.3 (Malliavin derivative's chain rule) :
 \begin{eqnarray*}
  \partial_x f_T(x,\sigma).v & = &
  \mathbb E[DF(X_{T}^{x}).D_{(h^{x,v},0)}^{\textrm{FV}}X_{T}^{x}]\\
  & = &
  \mathbb E[DF(X_{T}^{x}).\langle\mathbf DX_{T}^{x},I^{-1}(h^{x,v})\rangle_H]\\
  & = &
  \mathbb E[\langle\mathbf D(F\circ X_{T}^{x}),I^{-1}(h^{x,v})\rangle_H].
 \end{eqnarray*}
 \item Let $x\in\mathbb R^d$ be fixed. On one hand, for every $\varepsilon\in ]0,1]$ and $\sigma,\tilde\sigma\in\mathbb S_p$, by Taylor's formula :
 \begin{eqnarray*}
  \frac{\left|F(X_{T}^{\sigma +\varepsilon\tilde\sigma}) - F(X_{T}^{\sigma})\right|}{\varepsilon}
  & = &
  \left|\int_{0}^{1}
  \langle
  DF(X_{T}^{\sigma +\theta\varepsilon\tilde\sigma}),
  DX_{T}^{\sigma +\theta\varepsilon\tilde\sigma}.\tilde\sigma
  \rangle d\theta
  \right|\\
  & \leqslant &
  C_2\int_{0}^{1}
  \left(1 + \left\|X_{T}^{\sigma +\theta\varepsilon\tilde\sigma}\right\|\right)^N
  \left\|DX^{\sigma +\theta\varepsilon\tilde\sigma}.\tilde\sigma\right\|_{\infty;T}d\theta.
 \end{eqnarray*}
 At \cite{FV08}, Theorem 10.36, the constant involving in the upper-bound doesn't depend on the signal and on the collection of vector fields. At the second point of Corollary \ref{MAJsdrift}, the two constants involving in the upper-bound continuously depend on the $\gamma$-Lipschitz norm of the collection of vector fields. Then, there exists a constant $C_4 > 0$, depending on $\sigma$ and $\tilde\sigma$ but not on $\varepsilon$ and $\mathbb W$, such that for every $\theta\in [0,1]$,
 \begin{displaymath}
 \left\|DX^{\sigma +\theta\varepsilon\tilde\sigma}.\tilde\sigma\right\|_{\infty;T}
 \leqslant
 C_4e^{C_4M_{C_4,I,p}(\mathbb W)}
 \end{displaymath}
 and
 \begin{displaymath}
 \|X_{T}^{\sigma +\theta\varepsilon\tilde\sigma}\|
 \leqslant
 C_4\left(\|\mathbb W^g\|_{p\textrm{-var};T}\vee
 \|\mathbb W^g\|_{p\textrm{-var};T}^{p}\right).
 \end{displaymath}
 Since $W$ satisfies assumptions \ref{CSGRPgauss} and \ref{CMhyp}, by Proposition \ref{CLLrelation}, Theorem \ref{CLLprincipal}, the generalized Fernique's theorem (cf. Theorem \ref{GRPgauss}) and Cauchy-Schwarz's inequality :
 \begin{displaymath}
 \varepsilon\in ]0,1]
 \longmapsto
 \frac{\left|F(X_{T}^{\sigma +\varepsilon\tilde\sigma}) - F(X_{T}^{\sigma})\right|}{\varepsilon}
 \end{displaymath}
 is bounded by an integrable random variable not depending on $\varepsilon$. Therefore, by Lebesgue's theorem, $f_T(x,.)$ is differentiable on $\mathbb S_p$ and
 \begin{equation}\label{FDIFvega}
 \forall\sigma,\tilde\sigma\in\mathbb S_p
 \textrm{, }
 \partial_{\sigma} f_T(x,\sigma).\tilde\sigma =
 \mathbb E\left[\langle DF(X_{T}^{\sigma}),\partial_{\sigma} X_{T}^{\sigma}.\tilde\sigma\rangle\right].
 \end{equation}
 On the other hand, consider $\sigma,\tilde\sigma\in\mathbb S_p$ satisfying Assumption \ref{HYPsynth}. By construction, paths of the process $h^{\sigma,\tilde\sigma}$ are continuously differentiable from $[0,T]$ into $\mathbb R^d$ and $h_{0}^{\sigma,\tilde\sigma} = 0$. Then, $h^{\sigma,\tilde\sigma}$ is a $H^1$-valued random variable. By Duhamel's principle (cf. Theorem \ref{DFVsignal}) :
 \begin{eqnarray*}
  D_{(h^{\sigma,\tilde\sigma},0)}^{\textrm{FV}}X_{T}^{\sigma} & = &
  \int_{0}^{T}
  J_{T\leftarrow s}^{\mathbb W^g}\sigma\left(X_{s}^{\sigma}\right)dh_{s}^{\sigma,\tilde\sigma}\\
  & = &
  \partial_{\sigma} X_{T}^{\sigma}.\tilde\sigma.
 \end{eqnarray*}
 Therefore, by equality (\ref{FDIFvega}), Lemma \ref{H1DIFFrde} and \cite{NUALART06}, Proposition 1.2.3 (Malliavin derivative's chain rule) :
 \begin{eqnarray*}
  \partial_{\sigma} f_T(x,\sigma).\tilde\sigma & = &
  \mathbb E[DF(X_{T}^{\sigma}).D_{(h^{\sigma,\tilde\sigma},0)}^{\textrm{FV}}X_{T}^{\sigma}]\\
  & = &
  \mathbb E[DF(X_{T}^{\sigma}).\langle\mathbf DX_{T}^{\sigma},I^{-1}(h^{\sigma,\tilde\sigma})\rangle_H]\\
  & = &
  \mathbb E[\langle\mathbf D(F\circ X_{T}^{\sigma}),I^{-1}(h^{\sigma,\tilde\sigma})\rangle_H].
 \end{eqnarray*}
\end{enumerate}
\end{proof}
\noindent
\textbf{Remarks :}
\begin{enumerate}
 \item Note that the differentiability of partial maps of $f_T$ only requires assumptions \ref{CSGRPgauss}, \ref{CMhyp}, \ref{REGcoef} and \ref{Freg1}.
 \item Assumption \ref{CSGRPgauss} only guarantees the existence of the enhanced Gaussian process $\mathbb W$ at Theorem \ref{GRPgauss}. Therefore, when $p\in [1,2[$, Theorem \ref{CALCsensibilites} holds under the following assumption.
\end{enumerate}
%


%
\begin{assumption}\label{HYPMsynth}
The process $W$ is $d$-dimensional, Gaussian, centered and with continuous paths of finite $p$-variation ($p\in [1,2[$). Moreover
\begin{displaymath}
C_{0}^{1}\left([0,T];\mathbb R^d\right)
\subset H^1
\hookrightarrow
C^{p\textrm{-var}}\left([0,T];\mathbb R^d\right),
\end{displaymath}
and there exists a constant $C > 0$ such that :
\begin{displaymath}
\forall h\in C_{0}^{1}\left(
[0,T];\mathbb R^d\right)
\textrm{$,$ }
\|h\|_{H^1}\leqslant
C\|\dot h\|_{\infty;T}.
\end{displaymath}
Functions $\mu$ and $\sigma$ are $[p] + 2$ times differentiable, bounded and of bounded derivatives and, for every $a\in\mathbb R^d$, $\sigma(a)$ is an invertible matrix. Moreover, the function $\sigma^{-1} :\mathbb R^d\rightarrow\mathcal M_d(\mathbb R)$ is bounded.
\end{assumption}
\noindent
\textbf{Example.} A fractional Brownian motion $B^H$ of Hurst parameter $H\in ]1/2,1[$ satisfies Assumption \ref{HYPMsynth}. Indeed, by Kolmogorov's continuity criterion, there exists $\varepsilon > 0$ and $C_T\in L^1(\Omega;\mathbb R_{+}^{*})$ such that for every $(s,t)\in\bar\Delta_T$,
\begin{equation}\label{MBFholder}
\|B_{t}^{H} - B_{s}^{H}\|
\leqslant
C_T|t - s|^{1/p}
\textrm{ $\mathbb P$-a.s.}
\end{equation}
where $p := 1/(H - \varepsilon)\in [1,2[$.
%


%
\begin{corollary}\label{DIVsensibilites}
Under assumptions \ref{HYPMsynth} and, \ref{Freg1} or \ref{Freg2} :
\begin{enumerate}
 \item For every $x,v\in\mathbb R^d$, $I^{-1}(h^{x,v})\in\mathbb D^{1,2}(H)$ and
 \begin{displaymath}
 \partial_xf_T(x,\sigma).v =
 \mathbb E\left[
 F(X_{T}^{x})
 \delta\left[I^{-1}(h^{x,v})\right]\right].
 \end{displaymath}
 \item For every $x\in\mathbb R^d$ and every $\sigma,\tilde\sigma\in\mathbb S_{p + 1}$ satisfying Assumption \ref{HYPMsynth}, $I^{-1}(h^{\sigma,\tilde\sigma})\in\mathbb D^{1,2}(H)$ and
 \begin{displaymath}
 \partial_{\sigma}f_T(x,\sigma).\tilde\sigma =
 \mathbb E\left[
 F(X_{T}^{\sigma})
 \delta\left[I^{-1}(h^{\sigma,\tilde\sigma})\right]
 \right].
 \end{displaymath}
\end{enumerate}
\end{corollary}
%


%
\begin{proof}
Consider $x,v\in\mathbb R^d$ and $\sigma,\tilde\sigma\in\mathbb S_{p + 1}$ satisfying Assumption \ref{HYPMsynth}. If necessary, processes $h^{x,v}$ and $h^{\sigma,\tilde\sigma}$ introduced at Theorem \ref{CALCsensibilites} will be respectively denoted by $h^{x,v}(W)$ and $h^{\sigma,\tilde\sigma}(W)$.
\\
\\
With notations of Proposition \ref{DVITO2diff} :
\begin{displaymath}
h^{x,v} =
\frac{1}{T}
\int_{0}^{.}R_{s}^{1}(W,g)ds
\textrm{ and }
h^{\sigma,\tilde\sigma} =
\frac{1}{T}
\int_{0}^{.}R_{s}^{2}(W,g)ds
\end{displaymath}
where,
\begin{displaymath}
R^1(W,g) :=
\sigma^{-1}\left[\pi_V\left[0,x;(W,g)\right]\right]
D^1(W,g)v
\end{displaymath}
and
\begin{displaymath}
R^2(W,g) :=
\sigma^{-1}\left[\pi_V\left[0,x;(W,g)\right]\right]
D^1(W,g)D_{T}^{2}(W,g)D_{T}^{3}(W,g).
\end{displaymath}
First order derivatives of $R^1(w,g)$ and $R^2(w,g)$ with respect to the function $w\in C^{p\textrm{-var}}([0,T];\mathbb R^d)$, in the direction $\eta\in H^1\hookrightarrow C^{p\textrm{-var}}([0,T];\mathbb R^d)$, can be written as sums of products that only involve $D^1(w,g)$, $D^2(w,g)$ and $D^3(w,g)$, and their first order derivatives with respect to $w$, in the direction $\eta$.
\\
\\
Then, by Proposition \ref{DVITO2diff}, there exists two constants $C_1 > 0$ and $C_2 > 0$, not depending on $w$ and $g$, such that :
\begin{equation}\label{MAJdivs1}
\|\partial_w
R^1(w,g).\eta\|_{\infty;T}
\leqslant
C_1\exp\left[C_1(\|\eta\|_{H^1}^{p} +
\|w\|_{p\textrm{-var};T}^{p})\right]
\end{equation}
and
\begin{equation}\label{MAJdivs2}
\|\partial_w
R^2(w,g).\eta\|_{\infty;T}
\leqslant
C_2\exp\left[C_2(\|\eta\|_{H^1}^{p} +
\|w\|_{p\textrm{-var};T}^{p})\right],
\end{equation}
because
\begin{eqnarray*}
 \partial_w R^1(w,g).\eta & = &
 \partial_{(w,g)}
 R^1(w,g).(\eta,0)
 \textrm{ and}\\
 \partial_w R^2(w,g).\eta & = &
 \partial_{(w,g)}
 R^2(w,g).(\eta,0).
\end{eqnarray*}
On one hand, let show that $I^{-1}(h^{x,v})$ and $I^{-1}(h^{\sigma,\tilde\sigma})$ belong to $\mathbb D^{1,2}(H)$ :
\begin{enumerate}
 \item Firstly, for every $\varepsilon\in ]0,1]$ and $t\in [0,T]$, by Taylor's formula and inequality (\ref{MAJdivs1}) :
 \begin{eqnarray*}
  \frac{\left\|R_{t}^{1}[(W,g) + \varepsilon(\eta,0)] - R_{t}^{1}(W,g)\right\|}{\varepsilon}
  & = &
  \left\|
  \int_{0}^{1}
  DR_{t}^{1}\left[
  (W,g) +\varepsilon\theta(\eta,0)\right].(\eta,0)d\theta
  \right\|\\
  & \leqslant &
  C_3\exp\left[C_3(\|\eta\|_{H^1}^{p} +
  \|W\|_{p\textrm{-var};T}^{p})\right]
 \end{eqnarray*}
 where $C_3 > 0$ is a deterministic constant, not depending on $\varepsilon$, $t$, $\eta$ and $W$. Then, by Lebesgue's theorem :
 \begin{equation}\label{CDh1}
 D_{\eta}h^{x,v}(W) =
 \frac{1}{T}
 \int_{0}^{.}
 \partial_W R_{s}^{1}(W,g).\eta ds
 \end{equation}
 and $I^{-1}(h^{x,v})$ is locally derivable in the sense of Malliavin.
 \\
 \\
 Secondly, since $I$ is an isometry from $H$ into $H^1$, and $\sigma$ and $W$ satisfy Assumption \ref{HYPMsynth}, there exists two deterministic constants $C_4 > 0$ and $C_5 > 0$ such that :
 \begin{eqnarray*}
  \mathbb E\left[
  \|I^{-1}(h^{x,v})\|_{H}^{2}\right] & = &
  \mathbb E(
  \|h^{x,v}\|_{H^1}^{2})\\
  & \leqslant &
  C_4\mathbb E\left(\sup_{t\in [0,T]}\|\dot h_{t}^{x,v}\|^2\right)\\
  & \leqslant &
  C_5\mathbb E\left(\sup_{t\in [0,T]}\|J_{t\leftarrow 0}^{W^g}\|_{\mathcal M}^{2}\right).
 \end{eqnarray*}
 Then, by Corollary \ref{ITODPint}, $\mathbb E[\|I^{-1}(h^{x,v})\|_{H}^{2}] < \infty$.
 \\
 \\
 Thirdly, since $I : H\rightarrow H^1$ is an isometry, in particular $I^{-1}$ is a linear and continuously differentiable map :
 \begin{eqnarray*}
  \|\mathbf D[I^{-1}(h^{x,v})]\|_{H^{\otimes 2}}^{2} & \leqslant &
  \sup_{
  \begin{tiny}
  \begin{array}{rcl}
  a,b & \in & H\\
  \|a\otimes b\|_{H^{\otimes 2}} & \leqslant & 1
  \end{array}
  \end{tiny}}
  |\langle\mathbf D[I^{-1}(h^{x,v})],a\otimes b\rangle_{H^{\otimes 2}}|^2\\
  & = &
  \sup_{
  \begin{tiny}
  \begin{array}{rcl}
  a,b & \in & H\\
  \|a\otimes b\|_{H^{\otimes 2}} & \leqslant & 1
  \end{array}
  \end{tiny}}
  |\langle D_{I(a)}[I^{-1}(h^{x,v})](W),b\rangle_H|^2\\
  & = &
  \sup_{
  \begin{tiny}
  \begin{array}{rcl}
  a,b & \in & H\\
  \|a\otimes b\|_{H^{\otimes 2}} & \leqslant & 1
  \end{array}
  \end{tiny}}
  |\langle I^{-1}[D_{I(a)}h^{x,v}(W)],b\rangle_H|^2\\ 
  & = &
  \sup_{
  \begin{tiny}
  \begin{array}{rcl}
  a,b & \in & H\\
  \|a\otimes b\|_{H^{\otimes 2}} & \leqslant & 1
  \end{array}
  \end{tiny}}
  |\langle D_{I(a)}h^{x,v}(W),I(b)\rangle_{H^1}|^2\\
  & \leqslant &
  \sup_{
  \begin{tiny}
  \begin{array}{rcl}
  a & \in & H^1\\
  \|a\|_{H^1} & \leqslant & 1
  \end{array}
  \end{tiny}}
  \|D_ah^{x,v}(W)\|_{H^1}^{2}.
 \end{eqnarray*}
 Then, by equality (\ref{CDh1}), and since $\sigma$ and $W$ satisfy Assumption \ref{HYPMsynth} :
 \begin{eqnarray*}
  \mathbb E\left[\|\mathbf D[I^{-1}(h^{x,v})]\|_{H^{\otimes 2}}^{2}\right]
  & \leqslant &
  \frac{C_4}{T^2}\mathbb E\left(
  \sup_{
  \begin{tiny}
  \begin{array}{rcl}
  a & \in & H^1\\
  \|a\|_{H^1} & \leqslant & 1
  \end{array}
  \end{tiny}}
  \sup_{t\in [0,T]}
  \|\partial_WR_{t}^{1}(W,g).a\|^2\right)\\
  & \leqslant &
  \frac{C_{1}^{2}C_4}{T^2}\mathbb E\left[
  \exp\left[
  2C_1\left(1 +\|W\|_{p\textrm{-var};T}^{p}\right)\right]\right] < \infty
 \end{eqnarray*}
 by Fernique's theorem.
 \item Firstly, for every $\varepsilon\in ]0,1]$ and $t\in [0,T]$, by Taylor's formula and inequality (\ref{MAJdivs2}) :
 \begin{eqnarray*}
  \frac{\left\|R_{t}^{2}[(W,g) + \varepsilon(\eta,0)] - R_{t}^{2}(W,g)\right\|}{\varepsilon}
  & = &
  \left\|
  \int_{0}^{1}
  DR_{t}^{2}\left[
  (W,g) +\varepsilon\theta(\eta,0)\right].(\eta,0)d\theta
  \right\|\\
  & \leqslant &
  C_6\exp\left[C_6(\|\eta\|_{H^1}^{p} +
  \|W\|_{p\textrm{-var};T}^{p})\right]
 \end{eqnarray*}
 where $C_6 > 0$ is a deterministic constant, not depending on $\varepsilon$, $t$, $\eta$ and $W$. Then, by Lebesgue's theorem :
 \begin{equation}\label{CDh2}
 D_{\eta}h^{\sigma,\tilde\sigma}(W) =
 \frac{1}{T}
 \int_{0}^{.}
 \partial_W
 R_{s}^{2}(W,g).\eta ds
 \end{equation}
 and $I^{-1}(h^{\sigma,\tilde\sigma})$ is locally derivable in the sense of Malliavin.
 \\
 \\
 Secondly, since $I$ is an isometry from $H$ into $H^1$, and $\sigma$ and $W$ satisfy Assumption \ref{HYPMsynth}, there exists a deterministic constant $C_7 > 0$ such that :
 \begin{eqnarray*}
  \mathbb E\left[
  \|I^{-1}(h^{\sigma,\tilde\sigma})\|_{H}^{2}\right] & = &
  \mathbb E(
  \|h^{\sigma,\tilde\sigma}\|_{H^1}^{2})\\
  & \leqslant &
  C_4\mathbb E\left(\sup_{t\in [0,T]}\|\dot h_{t}^{\sigma,\tilde\sigma}\|^2\right)\\
  & \leqslant &
  C_7\mathbb E\left(\sup_{t\in [0,T]}
  \|D_{t}^{1}(W,g)\|_{\mathcal M}^{2}
  \|D_{T}^{2}(W,g)\|_{\mathcal M}^{2}
  \|D_{T}^{3}(W,g)\|^2\right).
 \end{eqnarray*}
 Then, by Corollary \ref{ITODPint}, $\mathbb E[\|I^{-1}(h^{\sigma,\tilde\sigma})\|_{H}^{2}] < \infty$.
 \\
 \\
 Thirdly, as at the first point, by equality (\ref{CDh2}), and since $\sigma$ and $W$ satisfy Assumption \ref{HYPMsynth} :
 \begin{eqnarray*}
  \mathbb E[\|\mathbf D[I^{-1}(h^{\sigma,\tilde\sigma})]\|_{H^{\otimes 2}}^{2}]
  & \leqslant &
  \frac{C_4}{T^2}\mathbb E\left(
  \sup_{
  \begin{tiny}
  \begin{array}{rcl}
  a & \in & H^1\\
  \|a\|_{H^1} & \leqslant & 1
  \end{array}
  \end{tiny}}
  \sup_{t\in [0,T]}
  \|\partial_WR_{t}^{2}(W,g).a\|^2\right)\\
  & \leqslant &
  \frac{C_{2}^{2}C_4}{T^2}\mathbb E\left[
  \exp\left[
  2C_2\left(1 +\|W\|_{p\textrm{-var};T}^{p}\right)\right]\right] < \infty
 \end{eqnarray*}
 by Fernique's theorem.
\end{enumerate}
Therefore, $I^{-1}(h^{x,v})$ and $I^{-1}(h^{\sigma,\tilde\sigma})$ belong to $\textrm{dom}(\delta)$ and, if $F$ satisfies Assumption \ref{Freg1}, by Theorem \ref{CALCsensibilites} :
\begin{displaymath}
\partial_xf_T(x,\sigma).v =
\mathbb E\left[
F(X_{T}^{x})
\delta\left[I^{-1}(h^{x,v})\right]\right]
\end{displaymath}
and
\begin{displaymath}
\partial_{\sigma}f_T(x,\sigma).\tilde\sigma =
\mathbb E\left[
F(X_{T}^{\sigma})
\delta\left[I^{-1}(h^{\sigma,\tilde\sigma})\right]
\right].
\end{displaymath}
On the other hand, consider a function $F :\mathbb R^d\rightarrow\mathbb R$ satisfying Assumption \ref{Freg2}, $(\varphi_n,n\in\mathbb N)$ a regularizing sequence of functions from $\mathbb R^d$ into $\mathbb R_+$ with compact supports, $K_n$ the support of $\varphi_n$, $F_n := \varphi_n\ast F$ and $f_{T}^{n}(x,\sigma) :=\mathbb E[F_n(X_T)]$ for every $n\in\mathbb N$.
\\
\\
Since $F$ satisfies Assumption \ref{Freg2}, there exists two constants $C_8 > 0$ and $N\in\mathbb N^*$, such that for every $a\in\mathbb R^d$,
\begin{equation}\label{CVDFreg2}
|F(a)|
\leqslant
C_8\left(1 + \|a\|\right)^N.
\end{equation}
So, $F\in L_{\textrm{loc}}^{1}(\mathbb R^d;\mathbb R)$ and the sequence of functions $(F_n,n\in\mathbb N)$ converges almost everywhere to $F$. In particular,
\begin{displaymath}
F_n\left(X_T\right)
\xrightarrow[n\rightarrow\infty]{\textrm{a.s.}}
F\left(X_T\right).
\end{displaymath}
By construction, for every $n\in\mathbb N$, $F_n$ satisfies Assumption \ref{Freg1}. Then,
\begin{displaymath}
\partial_xf_{T}^{n}(x,\sigma).v =
\mathbb E\left[
F_n(X_{T}^{x})
\delta\left[I^{-1}(h^{x,v})\right]\right]
\end{displaymath}
and
\begin{displaymath}
\partial_{\sigma}f_{T}^{n}(x,\sigma).\tilde\sigma =
\mathbb E\left[
F_n(X_{T}^{\sigma})
\delta\left[I^{-1}(h^{\sigma,\tilde\sigma})\right]
\right].
\end{displaymath}
By inequality (\ref{CVDFreg2}), for every $n\in\mathbb N$,
\begin{eqnarray*}
 \left|F_n(X_T)\right| & \leqslant &
 \int_{K_n}
 \varphi_n(a)\left|F(X_T - a)\right|da\\
 & \leqslant &
 \sup_{a\in K_n}(1 + \|X_T\| + \|a\|)^N
 \int_{K_n}\varphi_n(a)da\\
 & \leqslant &
 (1 + \|X_T\| + M)^N\in L^2(\Omega)
 \textrm{ with }
 M := \sup_{a\in K_n}\|a\|.
\end{eqnarray*}
Therefore, by Lebesgue's theorem :
\begin{displaymath}
f_{T}^{n}(x,\sigma)
\xrightarrow[n\rightarrow\infty]{}
f_T(x,\sigma).
\end{displaymath}
Moreover, since the random variables $\delta[I^{-1}(h^{x,v})]$ and $\delta[I^{-1}(h^{\sigma,\tilde\sigma})]$ belong to $L^2(\Omega)$, by Lebesgue's theorem :
\begin{displaymath}
\partial_xf_T(x,\sigma).v =
\mathbb E\left[
F(X_{T}^{x})
\delta\left[I^{-1}(h^{x,v})\right]\right]
\end{displaymath}
and
\begin{displaymath}
\partial_{\sigma}f_T(x,\sigma).\tilde\sigma =
\mathbb E\left[
F(X_{T}^{\sigma})
\delta\left[I^{-1}(h^{\sigma,\tilde\sigma})\right]
\right].
\end{displaymath}
\end{proof}
%


%
\section{Application to mathematical finance and simulations}
\noindent
In a first subsection, Theorem \ref{CALCsensibilites} and Corollary \ref{DIVsensibilites} are applied to the calculation of sensitivities in a financial market model with stochastic volatility, such that each equation is driven by a fractional Brownian motion of Hurst parameter belonging to $]1/4,1[$. In a second subsection, still with a fractional Brownian signal, simulations of the sensitivities with respect to the initial condition and to the collection of vector fields are provided, when the Hurst parameter of the fBm belongs to $]1/2,1[$.
%


%
\subsection{Computation of sensitivities in a fractional stochastic volatility model}
In this subsection, the prices process of risky assets is the solution of a fractional stochastic volatility model (taken in the sense of rough paths), and the sensitivity of an option's price to perturbations of the volatility is calculated by using Theorem \ref{CALCsensibilites} and Corollary \ref{DIVsensibilites}.
\\
\\
Consider a stochastic process $W$, and $\mu :\mathbb R^d\rightarrow\mathbb R^d$, $\kappa : \mathbb R^d\rightarrow\mathbb R_{+}^{d}$, $\sigma,\vartheta : \mathbb R^d\rightarrow\mathcal M_d(\mathbb R)$ and $F : \mathbb R^d\rightarrow\mathbb R_+$ five functions satisfying one of the two following assumptions.
%


%
\begin{assumption}\label{MVShyp1}
There exists two independent $d$-dimensional fBm $B^{H_1}$ and $B^{H_2}$, of respective Hurst parameters $H_1\in ]1/4,1[$ and $H_2\in ]1/4,1[$, such that $W := (B^{H_1},B^{H_2})$.
\\
\\
Functions $\mu$, $\sigma$ and $\vartheta$ satisfy Assumption \ref{HYPsynth} for $p := 1/(H_1 -\varepsilon)\vee1/(H_2 -\varepsilon) < 4$ and $\varepsilon > 0$ as close to $0$ as possible. Functions $\kappa$ and $F$ are such that $F\circ\kappa$ satisfies Assumption \ref{Freg1}.
\end{assumption}
%


%
\begin{assumption}\label{MVShyp2}
There exists two independent $d$-dimensional fBm $B^{H_1}$ and $B^{H_2}$, of respective Hurst parameters $H_1\in ]1/2,1[$ and $H_2\in ]1/2,1[$, such that $W := (B^{H_1},B^{H_2})$.
\\
\\
Functions $\mu$, $\sigma$ et $\vartheta$ satisfy Assumption \ref{HYPMsynth} for $p := 1/(H_1 -\varepsilon)\vee1/(H_2 -\varepsilon) < 2$ and $\varepsilon > 0$ as close to $0$ as possible. Functions $\kappa$ and $F$ are such that $F\circ\kappa$ satisfies Assumption \ref{Freg2}.
\end{assumption}
\noindent
Consider the financial market model consisting of $d$ risky assets, of prices $S_t$ at time $t\in [0,T]$ such that :
\begin{displaymath}
\left\{
\begin{array}{rcl}
S_t & := & \kappa\left(Y_t\right)\\
dY_t & = & \mu\left(Y_t\right)dt +\sigma\left(Z_t\right)dB_{t}^{H_1}\textrm{ ; }Y_0\in\mathbb R^d\\
dZ_t & = & \vartheta\left(Z_t\right)dB_{t}^{H_2}\textrm{ ; }Z_0\in\mathbb R^d
\end{array}
\right. ,
\end{displaymath}
and an option of payoff $F(S_T) := (F\circ\kappa)(Y_T)$ on these assets.
\\
\\
Consider $X := (Y,Z)$, $\mathbb W$ the enhanced Gaussian process associated to $W$, and $V := (V_1,\dots,V_{2d + 1})$ the collection of $\gamma$-Lipschitz vector fields on $\mathbb R_{1}^{d}\oplus\mathbb R_{2}^{d}$ ($\gamma > p$) defined by :
\begin{displaymath}
\forall a,b\in\mathbb R_{1}^{d}\oplus\mathbb R_{2}^{d}
\textrm{, }
\forall c\in\mathbb R
\textrm{, }
V(a)(b,c) :=
V_1(a)c +
V_2(a)b
\end{displaymath}
where,
\begin{displaymath}
V_1 :=
\begin{pmatrix}
\mu\circ\pi_{\mathbb R_{1}^{d}}\\
0
\end{pmatrix}
\textrm{ and }
V_2 :=
\begin{pmatrix}
\sigma\circ\pi_{\mathbb R_{2}^{d}} & 0\\
0 & \vartheta\circ\pi_{\mathbb R_{2}^{d}}
\end{pmatrix},
\end{displaymath}
and $\pi_{\mathbb R_{i}^{d}}$ is the canonical projection from $\mathbb R_{1}^{d}\oplus\mathbb R_{2}^{d}$ into $\mathbb R_{i}^{d}$ for $i = 1,2$.
\\
\\
Rigorously, $X = \pi_V(0,X_0;\mathbb W^g)$ with $\mathbb W^g := S_{[p]}(\mathbb W\oplus g)$ and $g :=\textrm{Id}_{[0,T]}$.
%


%
\begin{corollary}\label{MVSsensibilite}
With notations of Theorem \ref{CALCsensibilites} and Corollary \ref{DIVsensibilites} :
\begin{enumerate}
 \item Under Assumption \ref{MVShyp1}, $f_T$ is differentiable from $\mathbb S_p$ into $\mathbb R_+$ and, for every $\vartheta,\tilde\vartheta\in\mathbb S_p$ satisfying Assumption \ref{HYPsynth},
 \begin{displaymath}
 \partial_{\vartheta}f_T(\vartheta).\tilde\vartheta =
 \mathbb E\left[
 \langle\mathbf D(F\circ S_{T}^{\vartheta}) ;
 ((I_{H_1}^{-1}\circ\pi_{H_{1}^{1}})(h^{\vartheta,\tilde\vartheta}),
 (I_{H_2}^{-1}\circ\pi_{H_{2}^{1}})(h^{\vartheta,\tilde\vartheta}))\rangle_H
 \right]
 \end{displaymath}
 where for $i = 1,2$, $H_{i}^{1}$ is the Cameron-Martin's space of $B^{H_i}$, $\pi_{H_{i}^{1}}$ is the canonical projection from $H_{1}^{1}\oplus H_{2}^{1}$ into $H_{i}^{1}$, and
 \begin{displaymath}
 h^{\vartheta,\tilde\vartheta} :=
 \frac{1}{T}
 \int_{0}^{.}
 V_{2}^{-1}\left(X_{s}^{\vartheta}\right)
 J_{s\leftarrow T}^{\mathbb W^g}
 \partial_{\vartheta}X_{T}^{\vartheta}.\tilde\vartheta ds.
 \end{displaymath}
 \item Under Assumption \ref{MVShyp2}, for every $\vartheta,\tilde\vartheta\in\mathbb S_{p + 1}$ satisfying Assumption \ref{HYPMsynth}, $I^{-1}(h^{\vartheta,\tilde\vartheta})\in\mathbb D^{1,2}(H)$ and
 \begin{displaymath}
 \partial_{\vartheta}f_T(\vartheta).\tilde\vartheta =
 \mathbb E\left[
 F(S_{T}^{\vartheta})\left[
 \delta_{H_1}\left[(I_{H_1}^{-1}\circ\pi_{H_{1}^{1}})(h^{\vartheta,\tilde\vartheta})\right] +
 \delta_{H_2}\left[(I_{H_2}^{-1}\circ\pi_{H_{2}^{1}})(h^{\vartheta,\tilde\vartheta})\right]
 \right]\right].
 \end{displaymath}
\end{enumerate}
\end{corollary}
%


%
\begin{proof}
Straightforward application of Theorem \ref{CALCsensibilites} and of Corollary \ref{DIVsensibilites}.
\end{proof}
\noindent
\textbf{Remark.} Via Corollary \ref{IHexp} and its remark, it is possible to provide an explicit expression of the sensitivity of $f_T(\vartheta)$ with respect to $\vartheta$. Indeed, under Assumption \ref{MVShyp2}, for $i = 1,2$,
\begin{displaymath}
\delta_{H_i}\left[(I_{H_i}^{-1}\circ\pi_{H_{i}^{1}})(h^{\vartheta,\tilde\vartheta})\right] =
\delta_{1/2}\left[\left[
(\varphi_{H_i}D^{H_i - 1/2})\circ
(\varphi_{H_i}^{-1}D^1)\circ
\pi_{H_{i}^{1}}\right]
(h^{\vartheta,\tilde\vartheta})\right].
\end{displaymath}
%


%
\subsection{Simulations}
In order to simulate the sensitivities studied in this paper, let's first remind a result coming from \cite{LEJ10} on the convergence of the explicit Euler scheme associated to a differential equation driven by a $\alpha$-H\"older continuous function from $[0,T]$ into $\mathbb R^d$ ($\alpha\in ]1/2,1[$), taken in the sense of Young.
%


%
\begin{proposition}\label{LEJyoung}
Consider $x_0\in\mathbb R^d$, $w : [0,T]\rightarrow\mathbb R^d$ a $\alpha$-H\"older continuous function with  $\alpha\in ]1/2,1[$, and $V := (V_1,\dots,V_d)$ a collection of differentiable vector fields on $\mathbb R^d$ such that its derivative is $\gamma$-H\"older continuous from $\mathbb R^d$ into itself ($\gamma\in ]0,1[$ and $\gamma + 1 > 1/\alpha$). Then, there exists a constant $C > 0$ such that for every $n\in\mathbb N^*$,
\begin{displaymath}
\left\|x^n -\pi_V(0,x_0;w)\right\|_{\infty;T}
\leqslant
Cn^{1 - 2/p}
\end{displaymath}
where, $x^n$ is the step-$n$ explicit Euler scheme associated to $\pi_V(0,x_0;w)$ for the dissection $D^n := \{r_{k}^{n}\}\in D_T$ :
\begin{displaymath}
x_{t}^{n} :=
\sum_{k = 0}^{n - 1}
\left[
x_{k}^{n} +
\frac{x_{k + 1}^{n} - x_{k}^{n}}{r_{k + 1}^{n} - r_{k}^{n}}
(t - r_{k}^{n})
\right]\mathbf 1_{\left[r_{k}^{n},r_{k + 1}^{n}\right[}(t)
\textrm{ $;$ } t\in [0,T]
\end{displaymath}
with
\begin{displaymath}
\left\{
\begin{array}{rcl}
 x_{0}^{n} & := & x_0\\
 x_{k + 1}^{n} & = &
 x_{k}^{n} +
 V(x_{k}^{n})(w_{r_{k + 1}^{n}} - w_{r_{k}^{n}})
\end{array}
\right.
\end{displaymath}
for $k = 0,\dots,n - 1$.
\end{proposition}
\noindent
Refer to \cite{LEJ10}, Proposition 5 for a proof.
%


%
\begin{corollary}\label{LEJSyoung}
Consider $x_0\in\mathbb R$, a $1$-dimensional fractional Brownian motion $B^H$ of Hurst parameter $H\in ]1/2,1[$, $\mu$ and $\sigma$ two functions from $\mathbb R$ into $\mathbb R$ satisfying Assumption \ref{REGcoef} for $p := 1/(H -\varepsilon) < 2$ and $\varepsilon > 0$ as close to $0$ as possible, $V$ the vector field on $\mathbb R$ such that $V(a)(b,c) := \mu(a)c + \sigma(a)b$ for every $a,b,c\in\mathbb R$, $X := \pi_V(0,x_0;B^H)$, $Y :=\partial_xX^x.1$ and $Z := \partial_{\sigma}X^{\sigma}.\tilde\sigma$ for $\tilde\sigma\in\mathbb S_p$ arbitrarily chosen. Then, for every $r\geqslant 1$,
\begin{displaymath}
\lim_{n\rightarrow\infty}
\mathbb E\left(\left\|X^n - X\right\|_{\infty;T}^{r}\right) =
\lim_{n\rightarrow\infty}
\mathbb E\left(\left\|Y^n - Y\right\|_{\infty;T}^{r}\right) =
\lim_{n\rightarrow\infty}
\mathbb E\left(\left\|Z^n - Z\right\|_{\infty;T}^{r}\right) = 0
\end{displaymath}
where, for every $n\in\mathbb N^*$, $X^n$, $Y^n$ and $Z^n$ are respectively the explicit Euler schemes associated to $X$, $Y$ and $Z$ for the dissection $D^n := \{r_{k}^{n}\}\in D_T$. Moreover, the rate of convergence of each sequence is $n^{r(1-2/p)}$.
\end{corollary}
%


%
\begin{proof}
Processes $Y$ and $Z$ satisfy respectively :
\begin{displaymath}
Y =\pi_{A_1}(0,1;W^{\mu,\sigma})
\textrm{ and }
Z =\pi_{A_2}\left[0,0;(W^{\mu,\sigma},W^{\tilde\sigma})\right]
\end{displaymath}
where,
\begin{displaymath}
W^{\mu,\sigma} :=
\int_{0}^{.}\dot\mu\left(X_s\right)ds +
\int_{0}^{.}\dot\sigma\left(X_s\right)dB_{s}^{H}
\textrm{ and }
W^{\tilde\sigma} :=
\int_{0}^{.}
\tilde\sigma\left(X_t\right)dB_{t}^{H}
\end{displaymath}
and, $A_1$ and $A_2$ are two collections of affine vector fields on $\mathbb R$ defined by :
\begin{displaymath}
\forall a,b,c\in\mathbb R
\textrm{, }
A_1(a)b :=
ab
\textrm{ and }
A_2(a)(b,c) :=
ab + c.
\end{displaymath}
Since paths of $B^H$ are almost surely $1/p$-H\"older continuous by Kolmogorov's continuity criterion, \cite{FV08}, Theorem 6.8 implies that paths of $W^{\mu,\sigma}$ and $W^{\tilde\sigma}$ are also almost surely $1/p$-H\"older continuous. Then, $X$, $Y$ and $Z$ satisfy conditions of Proposition \ref{LEJyoung}, and there exists a random variable $C > 0$ such that for every $n\in\mathbb N^*$,
\begin{displaymath}
\left\|X^n -X\right\|_{\infty;T}
\textrm{, }
\left\|Y^n -Y\right\|_{\infty;T}
\textrm{ and }
\left\|Z^n -Z\right\|_{\infty;T}
\end{displaymath}
admit $Cn^{1 - 2/p}$ as upper-bound.
\\
\\
Finally, by reading carefully the proof of \cite{LEJ10}, Proposition 5, $C$ belongs to $L^r(\Omega)$ for every $r\geqslant 1$ by Fernique's theorem. Therefore, for every $r\geqslant 1$,
\begin{eqnarray*}
 \mathbb E(\|X^n -X\|_{\infty;T}^{r}) & \leqslant & \mathbb E(C^r)n^{r(1 - 2/p)}
 \xrightarrow[n\rightarrow\infty]{} 0\textrm{, }\\
 \mathbb E(\|Y^n -Y\|_{\infty;T}^{r}) & \leqslant & \mathbb E(C^r)n^{r(1 - 2/p)}
 \xrightarrow[n\rightarrow\infty]{} 0\textrm{ and }\\
 \mathbb E(\|Z^n -Z\|_{\infty;T}^{r}) & \leqslant & \mathbb E(C^r)n^{r(1 - 2/p)}
 \xrightarrow[n\rightarrow\infty]{} 0
\end{eqnarray*}
because $p < 2$.
\end{proof}
\noindent
\textbf{Remark.} About the approximation of the solution of SDEs driven by a fBm, refer also to A. Neuenkirch and I. Nourdin \cite{NN07}.
\\
\\
Let $n\in\mathbb N^*$ be fixed. With assumptions and notations of Corollary \ref{LEJSyoung}, at each iteration of step-$n$ explicit Euler schemes, the value of $B_{r_{. + 1}^{n}}^{H} - B_{r_{.}^{n}}^{H}$ is computed via the Wood-Chang's algorithm (cf. T. Dieker \cite{DIEKER04} about simulation methods of the fBm).
\\
\\
Let $F : \mathbb R\rightarrow\mathbb R$ be a function satisfying Assumption \ref{Freg1}. With notations of Section 3, in order to approximate $\partial_xf_T(x,\sigma).1$ (resp. $\partial_{\sigma}f_T(x,\sigma).\tilde\sigma$),
\begin{displaymath}
\mathbb E\left[\dot F\left(X_{T}^{n}\right)Y_{T}^{n}\right]
\textrm{ (resp. }
\mathbb E\left[\dot F\left(X_{T}^{n}\right)Z_{T}^{n}\right]
\textrm{)}
\end{displaymath}
is estimated by the empirical mean $\Theta_{m}^{n}(Y)$ (resp. $\Theta_{m}^{n}(Z)$) of the $m$-sample from the distribution of $F^Y := \dot F(X_{T}^{n})Y_{T}^{n}$ (resp. $F^Z := \dot F(X_{T}^{n})Z_{T}^{n}$). By Corollary \ref{LEJSyoung}, $F^Y$ and $F^Z$ belong to $L^2(\Omega)$. Then,
\begin{enumerate}
 \item By the strong law of large numbers :
 \begin{eqnarray*}
  \Theta_{m}^{n}(Y)
  & \xrightarrow[m\rightarrow\infty]{\textrm{a.s.}} &
  \theta^n(Y) :=
  \mathbb E\left[\dot F\left(X_{T}^{n}\right)Y_{T}^{n}\right]
  \textrm{ and}\\
  \Theta_{m}^{n}(Z)
  & \xrightarrow[m\rightarrow\infty]{\textrm{a.s.}} &
  \theta^n(Z) :=
  \mathbb E\left[\dot F\left(X_{T}^{n}\right)Z_{T}^{n}\right].
 \end{eqnarray*}
 \item By the central limit theorem and Slutsky's lemma :
 \begin{eqnarray*}
  \sqrt{n}
  \frac{\Theta_{m}^{n}(Y) - \theta^n(Y)}{s_{m}^{n}(Y)}
  & \xrightarrow[m\rightarrow\infty]{\mathcal D} &
  \mathcal N(0,1)
  \textrm{ and}\\
  \sqrt{n}
  \frac{\Theta_{m}^{n}(Z) - \theta^n(Z)}{s_{m}^{n}(Z)}
  & \xrightarrow[m\rightarrow\infty]{\mathcal D} &
  \mathcal N(0,1)
 \end{eqnarray*}
 where, $s_{m}^{n}(Y)$ (resp. $s_{m}^{n}(Z)$) is the empirical standard deviation of the $m$-sample from the distribution of $F^Y$ (resp. $F^Z$).
\end{enumerate}
At level $\alpha\in ]0,1[$, the second point provides the following confidence intervals :
\begin{displaymath}
\mathbb P\left[
\Theta_{m}^{n}(Y) - \frac{t_{\alpha}}{\sqrt{m}}s_{m}^{n}(Y)
\leqslant\theta^n(Y)\leqslant
\Theta_{m}^{n}(Y) + \frac{t_{\alpha}}{\sqrt{m}}s_{m}^{n}(Y)
\right]
\simeq
1 -\alpha
\end{displaymath}
and
\begin{displaymath}
\mathbb P\left[
\Theta_{m}^{n}(Z) - \frac{t_{\alpha}}{\sqrt{m}}s_{m}^{n}(Z)
\leqslant\theta^n(Z)\leqslant
\Theta_{m}^{n}(Z) + \frac{t_{\alpha}}{\sqrt{m}}s_{m}^{n}(Z)
\right]
\simeq
1 -\alpha
\end{displaymath}
where, $\Phi(t_{\alpha}) = 1 -\alpha /2$ and $\Phi$ is the distribution function of $\mathcal N(0,1)$.
\\
\\
\textbf{Example.} Assume that $T := 1$, $H := 0.6$, $n := 2^{15}$, $m := 500$, $\mu\equiv 0$, $\sigma : a\mapsto 1 + e^{-a^2}$, $\tilde\sigma : a\mapsto 1 +\pi/2 +\arctan(a)$, $F : a\mapsto a^2$ and $x := 1$ :
\begin{figure}[H]
\centering
\includegraphics[scale = 0.6]{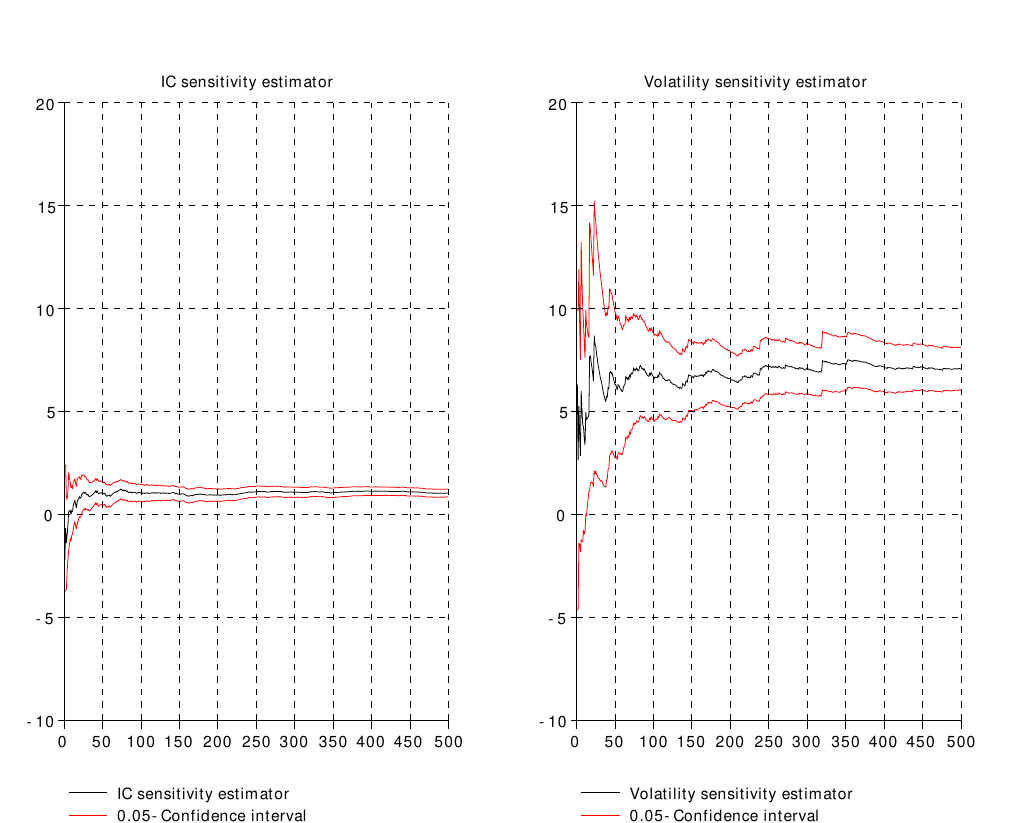}
\caption{Estimators convergence}
\end{figure}
\noindent
\begin{center}
\begin{tabular}{|l l|}
\hline
Statistics & Values\\
\hline\hline
$\Theta_{m}^{n}(Y)$ & $1.042$\\
$0.05$-confidence intervals & $[0.851;1.232]$\\
Length of the confidence interval & $0.381$\\
\hline
$\Theta_{m}^{n}(Z)$ & $7.112$\\
$0.05$-confidence intervals & $[6.071;8.154]$\\
Length of the confidence interval & $2.083$\\
\hline
\end{tabular}
\end{center}
%


%
\appendix
\section{Fractional Brownian motion}
\noindent
Essentially inspired by \cite{NUALART06} and \cite{DU98}, this appendix provides basics on the fractional Brownian motion, and the explicit expression of the isometry $I^{-1}$ defined at Lemma \ref{CMiso} for that Gaussian process.
%


%
\begin{definition}\label{fBm}
A fractional Brownian motion of Hurst parameter $H\in ]0,1]$ is a centered Gaussian process $B^H$ of covariance function $R_H$ defined by :
\begin{displaymath}
R_H(s,t) :=
\frac{1}{2}(
s^{2H} + t^{2H} - |t - s|^{2H})
\textrm{ $;$ }
s,t\in\mathbb [0,T].
\end{displaymath}
\end{definition}
\noindent
Let $B^H$ be a fractional Brownian motion of Hurst parameter $H\in ]0,1[$. The associated reproducing kernel Hilbert space is denoted by $\mathcal H$, the Wiener integral with respect to $B^H$ is denoted by $\mathbf B^H$, and the isometry provided at Lemma \ref{CMiso} is denoted by $I_H$.
%


%
\begin{definition}\label{OPfractionnaires}
Consider $\varphi :\mathbb R_+\rightarrow\mathbb R$ and $\alpha\in ]0,1]$ :
\begin{enumerate}
 \item If
 \begin{displaymath}
 l^{\alpha}(\varphi)(t) :=
 \frac{1}{\Gamma(\alpha)}\int_{0}^{t}
 (t - s)^{\alpha - 1}\varphi(s)ds
 \end{displaymath}
 exists for every $t\in\mathbb R_+$, $l^{\alpha}(\varphi)$ is the $\alpha$-fractional integral of $\varphi$.
 \item If
 \begin{displaymath}
 D^{\alpha}(\varphi)(t) :=
 \left\{
 \begin{array}{rcl}
  \displaystyle{\frac{1}{\Gamma(1 -\alpha)}\times\frac{d}{dt}\int_{0}^{t}(t - s)^{-\alpha}\varphi(s)ds}
  &\textrm{if}&
  \alpha\in ]0,1[\\
  \dot\varphi(t)
  &\textrm{if}&
  \alpha = 1
 \end{array}
 \right.
 \end{displaymath}
 exists for every $t\in\mathbb R_+$, $D^{\alpha}(\varphi)$ is the $\alpha$-fractional derivative of $\varphi$.
 \item If they are both defined :
 \begin{displaymath}
 (l^{\alpha}\circ D^{\alpha})(\varphi) =
 (D^{\alpha}\circ l^{\alpha})(\varphi) =
 \varphi.
 \end{displaymath}
\end{enumerate}
\end{definition}
\noindent
On fractional operators, refer to S. Samko et al. \cite{SKM93}.
\\
\\
\textbf{Notation.} $\mathcal E$ is the set of functions defined on $[0,T]$ by
\begin{displaymath}
\sum_{k = 1}^{n}
a_k\mathbf 1_{[0,s_k]}
\textrm{ ; }
n\in\mathbb N^*
\textrm{, }
(s_1,\dots,s_n)\in [0,T]^n
\textrm{, }
(a_1,\dots,a_n)\in\mathbb R^n.
\end{displaymath}
%


%
\begin{theorem}\label{MBFvolterra}
Let $K_{H}^{*}$ be the operator defined on $\mathcal E$ by :
\begin{displaymath}
\forall (s,t)\in\Delta_T
\textrm{$,$ }
K_{H}^{*}(\mathbf 1_{[0,t]})(s) :=
K_H(t,s)\mathbf 1_{[0,t]}(s)
\end{displaymath}
where,
\begin{displaymath}
K_H(t,s) :=
\frac{(t - s)^{H - 1/2}}{\Gamma(H + 1/2)}
\mathbf F\left(\frac{1}{2} - H,
H - \frac{1}{2},
H +\frac{1}{2},
1 -\frac{t}{s}\right)\mathbf 1_{[0,t[}(s)
\end{displaymath}
and $\mathbf F$ is the Gauss hyper-geometric function. Then,
\begin{enumerate}
 \item Let $J_H : L^2([0,T])\rightarrow H^1$ be the map defined by :
 \begin{displaymath}
 \forall h\in L^2([0,T])\textrm{$,$ }
 J_H(h) :=
 \int_{0}^{.}h(s)K_H(.,s)ds.
 \end{displaymath}
 For every $h\in\mathcal H$,
 \begin{displaymath}
 J_H(h) =
 \left\{
 \begin{array}{rcl}
 l^{2H}\circ (\varphi_{H}^{-1}l^{1/2 - H})\circ (\varphi_Hh) & \textrm{if} & H\leqslant 1/2\\
 l^1\circ (\varphi_Hl^{H - 1/2})\circ (\varphi_{H}^{-1}h) & \textrm{if} & H\geqslant 1/2
 \end{array}\right.
 \end{displaymath}
 where, $\varphi_H$ is the map defined by $\varphi_H(a) := a^{H - 1/2}\mathbf 1_{a\geqslant 0}$ for every $a\in\mathbb R$.
 \item The operator $K_{H}^{*}$ extends as an isometry from $\mathcal H$ into the closed subspace $\mathcal G := K_{H}^{*}(\mathcal H)$ of $L^2([0,T])$.
 \item The process $B := \mathbf B^H[(K_{H}^{*})^{-1}(\mathbf 1_{[0,.]})]$ is a standard Brownian motion, and
 \begin{displaymath}
 \forall t\in [0,T]\textrm{$,$ }
 B_{t}^{H} =
 \int_{0}^{t}
 K_H(t,s)dB_s.
 \end{displaymath}
 \item The divergence $\delta_H$ associated to $\mathbf B^H$ satisfies $\delta_H = \delta_{1/2}\circ K_{H}^{*}$.
 \end{enumerate}
\end{theorem}
\noindent
Refer to \cite{DU98}, Theorem 2.1 and Corollary 3.1, and \cite{NUALART06}, Proposition 5.2.2 for a proof.
\\
\\
\textbf{Remark.} At  \cite{DU98}, Theorem 3.3, L. Decreusefond and S. Ustunel proved that :
\begin{displaymath}
H^1 =
\left\{J_H(\dot h) ; \dot h\in L^2([0,T])\right\}.
\end{displaymath}
%


%
\begin{corollary}\label{IHexp}
The isometry $I_H$ satisfies $I_H = J_H\circ K_{H}^{*}$. In particular,
\begin{displaymath}
I_{H}^{-1} =
\left\{
\begin{array}{rcl}
(K_{H}^{*})^{-1}\circ (\varphi_{H}^{-1}D^{1/2 - H})\circ (\varphi_HD^{2H}) & \textrm{if} & H\leqslant 1/2\\
(K_{H}^{*})^{-1}\circ (\varphi_HD^{H - 1/2})\circ (\varphi_{H}^{-1}D^1) & \textrm{if} & H\geqslant 1/2
\end{array}\right. .
\end{displaymath}
\end{corollary}
%


%
\begin{proof}
On one hand, the isometry property of the It\^o's stochastic integral together with the third point of Theorem \ref{MBFvolterra} imply that for every $s,t\in [0,T]$,
\begin{displaymath}
\int_{0}^{s\wedge t}
K_H(t,u)K_H(s,u)du =
\mathbb E(B_{t}^{H}B_{s}^{H}).
\end{displaymath}
So, by definitions of $\mathbf B^H$ and $I_H$ :
\begin{displaymath}
\int_{0}^{s\wedge t}
K_{H}^{*}(\mathbf 1_{[0,t]})(u)K_H(s,u)du =
\mathbb E\left[\mathbf B^H(\mathbf 1_{[0,t]})B_{s}^{H}\right].
\end{displaymath}
Then, the construction of $I_H$ at Lemma \ref{CMiso} implies that :
\begin{displaymath}
I_H(\mathbf 1_{[0,t]}) =
(J_H\circ K_{H}^{*})(\mathbf 1_{[0,t]}).
\end{displaymath}
That equality extends on $\mathcal H$ by a classical continuity argument.
\\
\\
On the other hand, since $K_{H}^{*} : \mathcal H\rightarrow\mathcal G$ and $I_H : \mathcal H\rightarrow H^1$ are two invertible maps, the restriction $(J_H)_{|\mathcal G} = I_H\circ (K_{H}^{*})^{-1}$ is also invertible. Then, by the first point of Theorem \ref{MBFvolterra} :
\begin{displaymath}
I_{H}^{-1} =
\left\{
\begin{array}{rcl}
(K_{H}^{*})^{-1}\circ (\varphi_{H}^{-1}D^{1/2 - H})\circ (\varphi_HD^{2H}) & \textrm{if} & H\leqslant 1/2\\
(K_{H}^{*})^{-1}\circ (\varphi_HD^{H - 1/2})\circ (\varphi_{H}^{-1}D^1) & \textrm{if} & H\geqslant 1/2
\end{array}\right. .
\end{displaymath}
\end{proof}
\noindent
\textbf{Remark.} By the fourth point of Theorem \ref{MBFvolterra} and Corollary \ref{IHexp} :
\begin{eqnarray*}
 \delta_H\circ I_{H}^{-1} & = &
 \delta_{1/2}\circ (J_H)_{|\mathcal G}^{-1}.
\end{eqnarray*}
%


%

%
\end{document}